\documentclass[a4paper,11pt]{amsart}
\usepackage[english]{certus}

\makeatletter
\newcommand{\mylabel}[2]{#2\def\@currentlabel{#2}\label{#1}}
\makeatother

\DeclareMathOperator{\AFI}{AFI}
\DeclareMathOperator{\cyc}{cyc}
\DeclareMathOperator{\FinMod}{FMod}

\DeclareMathOperator{\HS}{HS}
\DeclareMathOperator{\hyp}{hyp}
\DeclareMathOperator{\LEF}{LEF}
\DeclareMathOperator{\PSL}{PSL}
\DeclareMathOperator{\sof}{sof}
\DeclareMathOperator{\Sub}{Sub}
\DeclareMathOperator{\weak}{weak}

\title{Sofic actions, halo products, and metric approximations of groups}
\author{Vadim Alekseev and Henry Bradford}

\begin{document}

\begin{abstract}
We introduce the notion of a ``sofic $\mathcal{C}$-action'' of one group 
on another by automorphisms, for $\mathcal{C}$ a class of groups. 
We show that if $\mathcal{C}$ is the class of (i) sofic, (ii) hyperlinear, 
(iii) linear sofic or (iv) weakly sofic groups, 
then the class $\mathcal{C}$ is closed under taking semidirect products 
with sofic $\mathcal{C}$-action. 
We use this to construct a wide variety of new examples of 
groups in the classes (i)-(iv), many of them arising as ``halo products'' 
in the sense of Genevois-Tessera. 
We have a parallel set of results producing new examples of semidirect 
products which are locally embeddable into finite groups. 
Our framework also unifies existing results in the literature, due to 
Hayes-Sale; Brude-Sasyk and Gao-Kunnawalkam Elayavalli-Patchell. 
\end{abstract}

\maketitle

\tableofcontents

\section{Introduction}

The search for a non-sofic group remains one of the major open problems in group theory. 
It even remains unknown whether the semidirect product $\Delta \rtimes_{\beta} \Gamma$ 
of two sofic groups $\Gamma$ and $\Delta$ is sofic. 
In this article we define the notion of a ``sofic action'' $\beta$ of $\Gamma$ on $\Delta$ 
by automorphisms, which is sufficient to guarantee soficity of 
$\Delta \rtimes_{\beta} \Gamma$ (provided $\Gamma$ and $\Delta$ are sofic). 
Our main technical result (Theorem \ref{BigTechThm}) yields the following. 

\begin{Thm} \label{BigTechThmIntro}
    Let $\mathcal{C}$ be the class of sofic groups. 
    Suppose that $\Gamma , \Delta \in \mathcal{C}$ and that 
    $\beta : \Gamma \rightarrow \Aut(\Delta)$ is a sofic $\mathcal{C}$-action by automorphisms. 
    Then $\Delta \rtimes_{\beta} \Gamma$ is a sofic group. 
\end{Thm}

See Definition \ref{SoficAutoActionDefn} for the notion of a sofic $\mathcal{C}$-action by automorphisms, for $\mathcal{C}$ a class of groups; 
our introduction of this concept is closely inspired by the notion 
of a sofic action of a group on a set \cite{GaoElayPatchSet} 
or a graph \cite{GaoElayPatchGraph}, recently introduced by Gao, 
Kunnawalkam Elayavalli and Patchell. 
We have corresponding results for sofic $\mathcal{C}$-actions, 
where $\mathcal{C}$ is the class of \emph{hyperlinear}; \emph{$\mathbb{K}$-linear sofic}, 
or \emph{weakly sofic} groups. 
Indeed our proof works for any class $\mathcal{C}$ of groups defined by the existence of metric approximations into a family $\mathcal{G}$ of groups equipped with bi-invariant metrics, 
such that $\mathcal{G}$ satisfies two assumptions which we call 
\emph{product-compatibility} and \emph{wreath-compatibility} (Definitions \ref{ProdCompatDef} and 
\ref{WrCompatDef}, respectively). 

Our main application of Theorem \ref{BigTechThmIntro} is to a class of groups called 
\emph{halo products}, recently introduced by Genevois and Tessera \cite{GeneTess}. 
The most familiar examples of halo products are (restricted) wreath products, 
but they also include \emph{free wreath products}; \emph{graph wreath products}; 
\emph{verbal wreath products}; \emph{symmetric enrichments} (of which the best known example is 
Houghton's second group), and many others (see Section 2 of \cite{GeneTess} 
for a large variety of examples). 
A halo product is constructed out of two ingredients. 
First, a \emph{halo of groups} $\mathcal{L}$ 
over a set, graph, or other structure $X$, 
which is an assignment, to every subset $Y \subseteq X$, a group $L(Y)$, 
subject to certain natural compatibility conditions. 
Second, an action $\alpha$ of a group $\Gamma$ on $X$ by automorphisms of the structure, 
inducing an action $\hat{\alpha}$ of $\Gamma$ by automorphisms on $L(X)$ 
which respects the halo of groups over $X$: that is, for $g \in \Gamma$ 
and $Y \subseteq X$ we have $\hat{\alpha} (L(Y)) = L(\alpha(g)(Y))$. 
We focus on haloes of groups over sets (respectively graphs) which satisfy a mild 
functoriality condition; we refer to such as \emph{set-theoretic} 
(respectively \emph{graph-theoretic}) haloes of groups. 

\begin{Prop} \label{SoficActSoficAutoPropIntro}
Let $\mathcal{C}$ be a class of groups 
and let $\mathcal{L}$ be a set-theoretic 
(respectively graph-theoretic) halo of groups. 
Suppose that for every finite set 
(graph) $X$ the group $L(X)$ lies in $\mathcal{C}$. 
Then for any set (graph) $X$ and and sofic group action 
$\Gamma \curvearrowright_{\alpha} X$, 
the induced action $\hat{\alpha} : \Gamma \rightarrow \Aut(L(X))$ 
is a sofic $\mathcal{C}$-action. 
\end{Prop}

The combination of Theorem \ref{BigTechThmIntro} and Proposition \ref{SoficActSoficAutoPropIntro} 
allow us to construct many families of sofic halo products, 
and by analogous methods we have many new examples of hyperlinear; linear sofic 
and weakly sofic groups; 
see Theorems \ref{SoficSymEnrThm}; \ref{LinEnrSoficThm}; 
\ref{GraphWPSoficThm}-\ref{GraphWPWSoficThm}; 
\ref{MainVerbalProdThm} and \ref{SoficAutEnrThm} below. 
These results generalise and unify several results in the literature on metric approximations of particular classes 
of halo products \cite{BrudSasy,ElekSzaboSofic,GaoElayPatchSet,GaoElayPatchGraph,HayesSaleMetric}. 
Here are two cases in which our conclusions are especially easy to state, 
for $\alpha$ an action of a group $\Gamma$ on a set. 

\begin{Ex}
For $\Delta$ a group, the \emph{(restricted permutational) wreath product} 
$\Delta \wr_X \Gamma$ of $\Delta$ and $\Gamma$ 
is the semidirect product $(\oplus_X \Delta) \rtimes \Gamma$, 
where $\Gamma$ acts on $\oplus_X \Delta$ by permuting co-ordinates 
according to $\alpha$. 
Hayes and Sale \cite{HayesSaleMetric} study the case of the \emph{regular} 
wreath product $\Delta \wr \Gamma$
(that for which $\alpha$ is the regular action of $\Gamma$ on itself)
and show that if $\Delta$ is sofic 
(respectively hyperlinear; $\mathbb{K}$-linear sofic; weakly sofic) 
then so is $\Delta \wr \Gamma$, provided $\Gamma$ is a sofic group. 
Meanwhile, Gao, Elayavalli and Patchell \cite{GaoElayPatchSet} 
show that if $\Gamma$ and $\Delta$ are sofic (respectively hyperlinear) 
then so is $\Delta \wr_X \Gamma$, whenever $\alpha$ is a sofic action. 
Our results extend those of \cite{GaoElayPatchSet}, 
showing in addition that if $\Gamma$ and $\Delta$ are $\mathbb{K}$-linear sofic
(respectively weakly sofic) then so is $\Delta \wr_X \Gamma$, 
whenever $\alpha$ is a sofic action. 
\end{Ex}

\begin{Ex}
Let $\Sym_f (X)$ be the group of permutations of $X$ with finite support. 
Then $\alpha$ induces an action of $\Gamma$ on $\Sym_f (X)$ by automorphisms; 
the \emph{symmetric enrichment} of the action $\Gamma \curvearrowright_{\alpha} X$ 
is by definition the associated semidirect product $\Sym_f (X) \rtimes \Gamma$. 
Our results show that if $\Gamma$ is sofic 
(respectively hyperlinear; $\mathbb{K}$-linear sofic; weakly sofic), 
then so is $\Sym_f (X) \rtimes \Gamma$, for every sofic action $\alpha$. 
    Here our result is new even for the case of a regular action of a sofic group $\Gamma$ 
    on itself. 
\end{Ex}

Concurrently, we introduce sufficient criteria for a semidirect product 
to be \emph{locally embeddable into finite groups} (or \emph{LEF} for short). 
To this end we introduce a notion of ``LEF action'' $\beta$ 
of a group $\Gamma$ by automorphisms on a group $\Delta$ (see Definition \ref{LEFActGrpDef} below), 
generalising Vershik and Gordon's notion 
of an \emph{equivariantly approximable} action \cite{VersGord}. 
Parallel to Theorem \ref{BigTechThmIntro} we prove the following (Theorem \ref{VerGorEAProp}). 

\begin{Thm} \label{VerGorEAPropIntro}
    Let $\mathcal{C}$ be the class of LEF groups; let $\Gamma$ be a LEF group and suppose 
$\Gamma \curvearrowright_{\beta} \Delta$ is an $\mathcal{C}$-LEF action by automorphisms. 
Then $\Delta \rtimes_{\beta} \Gamma$ is a LEF group. 
\end{Thm}

Once again, our main application of Theorem \ref{VerGorEAPropIntro} is to halo products. 
A notion of a LEF action of a group on a set/graph was introduced in \cite{BradLEFExtn} 
(see Definition \ref{LEFActionDef} below). 
In fact, a transitive action of a group is LEF iff a point-stabiliser 
is a subgroup ``almost of finite index'' in the sense of Cornulier \cite{CornulierWreath} 
(see Proposition \ref{AFILEFActProp} below). 
Parallel to Proposition \ref{SoficActSoficAutoPropIntro} we have the following 
(Proposition \ref{LEFActEAAutoProp}). 

\begin{Prop} \label{LEFActEAAutoPropIntro}
Let $\mathcal{L}$ be a set-theoretic
(respectively graph-theoretic) halo of groups and let $\mathcal{C}$ be a class of groups. 
Suppose that for every finite set (graph) $Y$, $L(Y)$ lies in $\mathcal{C}$. 
Let $X$ be a set (graph) and suppose $\Gamma \curvearrowright_{\alpha} X$ is a LEF action. 
Then the induced action $\hat{\alpha} : \Gamma \rightarrow \Aut(L(X))$
is a $\mathcal{C}$-LEF action. 
\end{Prop}

Along the way, we prove the following (Theorem \ref{LERFLEFActThm}), which strengthens a recent result of Gao \cite{Gao}. 

\begin{Thm} 
    If $\Gamma$ is LERF, then all actions of $\Gamma$ on sets are LEF actions. 
\end{Thm}

Combining Theorem \ref{VerGorEAPropIntro} and Proposition \ref{LEFActEAAutoPropIntro} we 
obtain many families of examples of halo products which are LEF 
(see Theorems \ref{SymEnrLEFThm}; \ref{LinEnrLEFThm}; \ref{GraphWPLEFThm}; 
\ref{LEFVerbalProdThm} and \ref{LEFAutEnrThm} below). 

All the approximation properties discussed hitherto (LEF, soficity, linear soficity, hyperlinearity, weak soficity) are proper generalisations of residual finiteness. Where known, we also comment on the conditions under which residual finiteness holds for the examples of halo products we consider in Section \ref{ExamplesSect}. 
As a rule, these conditions are very restrictive: 
see Proposition \ref{SymEnrRFProp} (for symmetric and alternating enrichments); 
Proposition \ref{ElemNoRFProp} (for elementary enrichments), 
and Theorem \ref{NeedhamGWPRFThm} (for graph wreath products, due to Needham \cite{Needham}).

\addtocontents{toc}{\protect\setcounter{tocdepth}{0}}
\subsection*{Notation}

For $\Gamma$ a group, we write $\Delta \leq_f \Gamma$ to denote a finite-index subgroup $\Delta$ of $\Gamma$. 
For $S$ a subset of a group $\Gamma$ we denote by $S^{-1}$ the set of inverses of the elements 
of $S$. 
For $l \in \mathbb{N}$, 
we denote by $B_S(l)$ the set of all elements of $\Gamma$ 
expressible as a product, of length at most $l$, 
of the elements of $S$ and their inverses. 

\addtocontents{toc}{\protect\setcounter{tocdepth}{2}}
\section{Halo products}

\begin{Def} \label{HaloDefn}
A \emph{halo of groups} over a set $X$ 
is a family of groups $\mathcal{L}_X = \lbrace L(Y) : Y \subseteq X \rbrace$ 
satisfying the following axioms: 
\begin{itemize}
\item[(i)] For all $Z \subseteq Y \subseteq X$ then $L(Z) \leq L(Y)$; 

\item[(ii)] $L(\emptyset)$ is the trivial group, 
and $L(X) = \langle L(Y) : Y \subseteq X \text{ finite} \rangle$; 

\item[(iii)] For all $Y,Z \subseteq X$, $L(Y \cap Z) = L(Y) \cap L(Z)$. 

\end{itemize}
\end{Def}

\begin{Rem}
Is is an easy consequence of axioms (i)-(iii) that for any $Y \subseteq X$, 
$L(Y) = \langle L(Z) : Z \subseteq Y \text{ finite} \rangle$. 
\end{Rem}

\begin{Def} \label{HaloProdDefn}
Given a group action $\Gamma \curvearrowright_{\alpha} X$ and 
a homomorphism $\hat{\alpha} : \Gamma \rightarrow \Aut(L(X))$ satisfying: 
$$\hat{\alpha}(g) \big( L(Y) \big) = L(\alpha(g)(Y)) 
\text{ for all }g \in \Gamma \text{ and } Y \subseteq X$$
we refer to the semidirect product 
$L(X) \rtimes_{\hat{\alpha}} \Gamma$ as the associated 
\emph{halo product} of groups. 
\end{Def}

We refer the reader to \cite{GeneTess} Section 2 for numerous 
natural examples of haloes of groups 
and associated halo products, 
several of which shall be revisited in Section \ref{ExamplesSect}. 
In many of the examples explored in \cite{GeneTess}, 
$X$ has the structure of a set; a graph, or a totally ordered set, 
and the group-theoretic structure of $L(Y) \leq L(X)$ 
for $Y \subseteq X$ depends only on the cardinality of $Y$; 
the structure of $Y$ as an induced subgraph of $X$, 
or the induced order on $Y$, respectively, 
in such a manner that \emph{any} 
group action $\Gamma \curvearrowright_{\alpha} X$ 
(preserving the relevant structure) induces an 
action of $\Gamma$ on $L(X)$ by automorphisms satisfying the 
conditions of Definition \ref{HaloProdDefn}. 
We formalise this idea as follows. 

\begin{Def} \label{FunctHaloDefn}
A \emph{set-theoretic} (respectively \emph{graph-theoretic}; 
\emph{digraph-theoretic}) 
halo of groups is a covariant functor $\mathcal{L}$ from the category 
of sets equiped with injective functions 
(from the category of graphs 
equipped with inclusions of induced 
subgraphs; the category of digraphs equipped with 
inclusions of induced subgraphs) 
to the category of groups equipped with monomorphisms, 
such that for any set (graph; digraph) $X$, 
the family $\lbrace \hat{\iota}(Y) : Y \subseteq X \rbrace$ 
forms a halo of groups over $X$. 
Here for $\pi : Y \hookrightarrow X$ an injective map in the 
relevant category, we write $\hat{\pi} : L(Y) \hookrightarrow L(X)$ 
for the monomorphism of groups induced by the functor, 
and for $Y \subseteq X$, $\iota : Y \hookrightarrow X$ 
is the inclusion map. 
\end{Def}

\begin{Def} \label{FunctHaloProdDefn}
If $\mathcal{L}$ is a set-theoretic (respectively graph-theoretic; 
digraph-theoretic) halo of groups, 
then for any set (graph; digraph) $X$, 
and any group action $\Gamma \curvearrowright_{\alpha} X$ 
(preserving the relevant structure), 
there is an induced action of $\Gamma$ on $L(X)$ by automorphisms 
given by: 
$$g \mapsto (\tau \mapsto \hat{\alpha}_g (\tau)) \text{ (for $g \in \Gamma$ 
and $\tau \in L(X)$)}$$
such that the conditions of Definition \ref{HaloProdDefn} are satisfied. 
We write $\mathcal{L}(\Gamma,X) = L(X) \rtimes_{\hat{\alpha}} \Gamma$ 
for the associated halo product. 
\end{Def}

\begin{Rem} \label{InclusionOrbitRmrk}
If $Y \subseteq X$ is a $\Gamma$-invariant subset (induced subgraph), 
then there is an obvious inclusion $\mathcal{L}(\Gamma,Y)$ 
as a subgroup of $\mathcal{L}(\Gamma,X)$. 
\end{Rem}

\begin{Ex}
Here is an example of a halo of groups discussed in \cite{GeneTess} Section 2, 
which is built out of graphs but is not 
a graph-theoretic halo in the sense of Definition \ref{FunctHaloDefn}. 
To a graph $X=(V,E)$ one can associate a surface $\Sigma_X$, 
which is obtained as a connected sum of tori indexed by $V$, 
with the tubes between the tori being pairwise disjoint and indexed by $E$ 
(that is, for each $e=(v,w) \in E$ we join the tori corresponding to $v,w \in V$ 
by a tube). 
Let $\FinMod(\Sigma_X)$ be the compactly-supported mapping class group of $\Sigma_X$, 
and for $Y = (V',E')$ an induced subgraph of $X$, 
let $L(Y) \leq \FinMod(\Sigma_X)$ denote the subgroup consisting of mapping classes 
realisable by homeomorphisms of $\Sigma_X$ supported on the tori 
corresponding to $V'$ and the tubes corresponding to $E'$ 
(that is, they fix pointwise any tori corresponding to $V \setminus V'$ 
and any tubes corresponding to $E \setminus E'$). 
Then $L(X) =  \FinMod(\Sigma_X)$ and 
$\mathcal{L} = \lbrace L(Y) : Y \subseteq X \rbrace$ is a halo of groups over $X$. 
But this is not a graph-theoretic halo: 
if it were, then for $X$ the path of length three and $Z$ the path of length two, 
the inclusion of $Z$ into $X$ as an induced subgraph would induce a monomorphism of 
groups $\FinMod (\Sigma_Z) \hookrightarrow \FinMod (\Sigma_X)$. 
But $\Sigma_X = \Sigma_4$ and $\Sigma_Z = \Sigma_3$ are the closed orientable surfaces 
of genus four and three, respectively, 
so $\Mod (\Sigma_3) \hookrightarrow \Mod (\Sigma_4)$, 
contradicting Theorem 1 of \cite{ChenLani}. 
\end{Ex}

\section{Metric approximations of groups}

\subsection{General framework} \label{GenMetricApproxSubsect}

Let $\Gamma$ be a group and let $G$ be a group
equipped with a bounded bi-invariant metric $d$. 
Let $\varphi : \Gamma \rightarrow G$ be a function,
let $F \subseteq \Gamma$ be a finite subset and let $\epsilon > 0$.

\begin{itemize}
\item[(i)] We describe $\varphi$ as \emph{$(F,\epsilon)$-multiplicative} if,
for all $g,h \in F$ such that $gh \in F$,
$$d\big( \varphi(gh),\varphi(g)\varphi(h) \big) < \epsilon; $$
\item[(ii)] For $c>0$,
an element $g \in G$ is \emph{$c$-separated}
if $d(g,e_G) > c$.
We describe $\varphi$ as \emph{$(F,c)$-separating} if,
for all $e \neq g \in F$, $\varphi(g) \in G$ is $c$-separated;
\item[(iii)] $\varphi$ is \emph{unital} if 
$\varphi(e_{\Gamma}) = e_G$.
\end{itemize}

\begin{Def}
A \emph{metric family of groups} is a family
$\mathcal{G} = \lbrace (G_i,d_i) : i \in I \rbrace$,
where for all $i \in I$, $G_i$ is a nontrivial group and
$d_i$ is a bi-invariant metric on $G_i$, 
such that the diameters of the metric spaces $(G_i,d_i)$ are uniformly bounded. 
If $F \subseteq \Gamma$, $\epsilon,c>0$ 
are such that the function $\varphi : \Gamma \rightarrow G_i$
(for some $i \in I$) is $(F,\epsilon)$-multiplicative;
$(F,c)$-separating and unital, 
we shall call it an \emph{$(F,\epsilon,c)$-$\mathcal{G}$-representation of
$\Gamma$}.
The group $\Gamma$ is a \emph{$(\mathcal{G},c)$-approximable group} if 
for all $F \subseteq \Gamma$ finite and $\epsilon > 0$
there exists an 
$(F,\epsilon,c)$-$\mathcal{G}$-representation of $\Gamma$. 
We refer to $\Gamma$ as \emph{$\mathcal{G}$-approximable} 
if it is $(\mathcal{G},c)$-approximable for some $c>0$. 
\end{Def}

\begin{Lemma} \label{BoringUnitalLemma}
Let $\Gamma$ be a group and $(G,d)$ be a metric group. 
Let $e_{\Gamma} \in F \subseteq\Gamma$; $\delta , c > 0$. 
Suppose $\varphi : \Gamma \rightarrow G$ satisfies: 
\begin{itemize}
    \item[(i)] $\varphi$ is $(F,\delta)$-multiplicative; 

    \item[(ii)] $d(\varphi(e_{\Gamma}),e_G) < \delta$; 

    \item[(iii)] For all $g,h \in F$ with $g \neq h$, $d(\varphi(g),\varphi(h)) \geq c$. 
    
\end{itemize}
Then there exists an $(F,2\delta,c-\epsilon)$-$\mathcal{G}$-representation 
$\varphi ' : \Gamma \rightarrow G$. 
\end{Lemma}

\begin{proof}
Define $\varphi '$ by $\varphi ' (e_{\Gamma}) = e_G$ and $\varphi '(g) = \varphi(g)$ 
for $g \neq e_{\Gamma}$. Then $\varphi '$ is unital. 
We check $(F,2\delta)$-multiplicativity; let $g,h \in F$ be such that $gh \in F$. 
If $g,h,gh \neq e_{\Gamma}$ there is nothing to check. 
For $g = e_{\Gamma}$ or $h = e_{\Gamma}$ we have 
$\varphi '(g)\varphi '(h) = \varphi '(gh)$. 
In the remaining case, $h = g^{-1} \neq e$, so that: 
$$d(\varphi '(g)\varphi '(h),\varphi '(gh)) = d(\varphi (g)\varphi (g^{-1}),e_G)
\leq d(\varphi (g)\varphi (g^{-1}),\varphi (e_{\Gamma})) + d(\varphi (g)\varphi (\varphi (e_{\Gamma}),e_G) < 2 \delta. $$
Finally, for $e \neq g \in F$, 
$$d(\varphi '(g),e_G) \geq d(\varphi(g),\varphi (e_{\Gamma})) - d(\varphi (e_{\Gamma}),e_G) 
> c-\delta$$
so $\varphi'$ is $(F,c-\delta)$-separating. 
\end{proof}

\begin{Rem} \label{LocalMetApproxRem}
    Since the existence of an $(F,\epsilon,c)$-$\mathcal{G}$-representation of
$\Gamma$ depends only on the finite subset $F$ of $\Gamma$, 
the group $\Gamma$ is $\mathcal{G}$-approximable iff every finitely generated 
subgroup of $\Gamma$ is $\mathcal{G}$-approximable. 
\end{Rem}

The next observation is stated as Lemma 2.8 of \cite{HayesSaleMetric} 
in the special case of wreath products, but holds more generally. 

\begin{Lemma} \label{GenHayesSaleLem}
Let $H = K\rtimes_{\alpha} \Gamma$ be a semidirect product of two discrete groups. For every finite set $F\subset K\rtimes \Gamma$ there are finite sets $F_1\subset K$ and $F_2\subset \Gamma$ such that the following holds. Let $\eps > 0$ and let $G$ be a group equipped with a bounded 
bi-invariant metric $d$. Suppose $\sigma\colon K\rtimes_{\alpha} \Gamma \to G$ be a map such that:
\begin{itemize}
\item[(i)] the restriction of $\sigma$ to $K$ is $(F_1,\eps/6)$-multiplicative;
\item[(ii)] the restriction of $\sigma$ to $\Gamma$ is $(F_2,\eps/6)$-multiplicative;
\item[(iii)] $\max_{k\in F_1,g\in F_2} d(\sigma(k,g),\sigma(k,1)\sigma(1,g)) < \eps/6$;
\item[(iv)]  $\max_{k\in F_1,g\in F_2} d(\sigma(1,g)\sigma(k,1),\sigma(\alpha(g)[k],1)\sigma(1,g)) < \eps/6$.
\end{itemize}
Then $\sigma$ is $(F,\eps)$-multiplicative.
\end{Lemma}

\begin{proof}
    Let $\lvert F \rvert = n$ and let $g_i \in \Gamma$, $k_i \in K$ be such that 
    $F = \lbrace (k_1,g_1),\ldots,(k_n,g_n) \rbrace$. 
    Set $F_2 = \lbrace e,g_1 ,\ldots,g_n\rbrace$ and set 
    $F_1 = \lbrace \alpha(g)[k_i] : 1 \leq i \leq n , g \in F_2 \rbrace$. 
    For $1 \leq i,j \leq n$ we are comparing $\sigma(k_i,g_i)\sigma(k_j,g_j)$ with 
    $\sigma((k_i,g_i)(k_j,g_j)) = \sigma(k_i \alpha(g_i)[k_j],g_i g_j)$ 
    and wish to show that they are at distance $< \epsilon$ in $G$. 
    This is immediate from repeated application of the triangle inequality. 
\end{proof}

\begin{Rem} \label{GenHayesSaleBiggerRmrk}
Note that conditions (i)-(iv) of Lemma \ref{GenHayesSaleLem}
only become stronger upon enlarging $F_1$ and $F_2$.
Therefore, if $F_1 \subseteq K$ and $F_2 \subseteq \Gamma$
satisfy the conclusion of Lemma Lemma \ref{GenHayesSaleLem} 
with respect to the finite set $F \subseteq K \rtimes \Gamma$,
then so too do any finite sets $F_1 ' \subseteq K$ and $F_2 ' \subseteq \Gamma$
satisfying $F_1 \subseteq F_1 '$ and $F_2 \subseteq F_2 '$.
\end{Rem}

\begin{Def} \label{ProdCompatDef}
The metric family $\mathcal{G}$ is \emph{product-compatible} 
if there exists $\mu>0$ such that for all $i,j \in I$, there exists $k \in I$
and a function $\Delta_{i,j} : G_i \times G_j \rightarrow G_k$ 
such that 
\begin{itemize}
\item[(a)] For all $\epsilon > 0$ there exists $\delta > 0$ 
such that for all $i,j \in I$, all $g_i , g_i ' \in G_i$
and all $g_j , g_j ' \in G_j$, 
if $d_i (g_i , g_i ') , d_j (g_j , g_j ') < \delta$
then $d_k (\Delta_{i,j}(g_i,g_j),\Delta_{i,j}(g_i ',g_j ')) < \epsilon$;

\item[(b)] For all $i,j \in I$, all $g_i , g_i ' \in G_i$
and all $g_j , g_j ' \in G_j$, $$d_k (\Delta_{i,j}(g_i,g_j),\Delta_{i,j}(g_i ',g_j ')) 
\geq \mu \max \lbrace d_i (g_i , g_i ') , d_j (g_j , g_j ') \rbrace.$$
\end{itemize}
We shall call $\Delta_{i,j}$ a \emph{$\mu$-product map} 
(or simply a product map if the specific value of $\mu$ is irrelevant). 
\end{Def}

\begin{Lemma} \label{ProdMapAlmostHomLem}
Suppose $\mathcal{G}$ is a product-compatible metric family, 
and suppose $\Delta_{i,j} : G_i \times G_j \rightarrow G_k$ is a product map. 
For all $\epsilon > 0$ there exists $\delta > 0$ such that, 
for any group $\Gamma$; $F \subseteq \Gamma$ and functions 
$\phi_i : \Gamma \rightarrow G_i$, 
$\phi_j : \Gamma \rightarrow G_j$, 
if $\phi_i , \phi_j$ are $(F,\delta)$-multiplicative, 
then $\Delta_{i,j} \circ (\phi_i \times \phi_j) : \Gamma \rightarrow G_k$ 
is $(F,\epsilon)$-multiplicative. 
\end{Lemma}

\begin{proof}
Given $\epsilon > 0$, let $\delta > 0$ be as in Definition \ref{ProdCompatDef} (a). 
    Let $g,h \in F$ with $gh \in F$, so that: 
$$d_i (\phi_i (gh),\phi_i(g)\phi_i(h)) , d_j (\phi_j (gh),\phi_j(g)\phi_j(h)) < \delta.$$
Now since $\Delta_{i,j}$ is a homomorphism, 
we have $$\Delta_{i,j} (\phi_i(g)\phi_i(h),\phi_j(g)\phi_j(h))
 = \Delta_{i,j} (\phi_i(g),\phi_j(g))  \Delta_{i,j} (\phi_i(h),\phi_j(h))$$
 so that by Definition \ref{ProdCompatDef} (a), 
 $$d_k (\Delta_{i,j}(\phi_i (gh),\phi_j (gh))  ,\Delta_{i,j} (\phi_i(g),\phi_j(g))  \Delta_{i,j} (\phi_i(h),\phi_j(h)) ) < \epsilon$$
 as desired. 
\end{proof}

\begin{Def} \label{WrCompatDef}
The metric family $\mathcal{G}$ is \emph{wreath-compatible}
if for all $i \in I$ and $n \in \mathbb{N}$ 
there exists $j \in J$ and homomorphisms 
$\tau : G_i ^n \rightarrow G_j$, $\psi : \Sym(n) \rightarrow G_j$
such that for all $\epsilon > 0$ there exists $\delta > 0$ 
such that for all $c > 0$, 
all $g_x , g_x ' \in G_i$ (with $x \in [n]$)
and $\sigma_1 , \sigma_2 \in \Sym(n)$,
\begin{itemize}
\item[(a)] If there is a subset $X \subseteq [n]$ such that
$\lvert X \rvert > (1-\delta) n$
and for all $x \in X$, $d_i (g_x , g_x ') < \delta$
(respectively $d_i (g_x , g_x ') > c$)
then: 
\begin{center}
    $d_j (\tau((g_x)_x),\sigma((g_x ')_x)) < \epsilon$
(respectively $d_j (\sigma((g_x)_x),\sigma((g_x ')_x)) > (1 - \epsilon) c$);
\end{center}

\item[(b)] If $d_n (\sigma_1 , \sigma_2) < \delta$
where $d_n$ is the normalized Hamming metric on $\Sym(n)$ 
(see Definition \ref{HammingDef} below), 
then $d_j (\psi(\sigma_1),\psi(\sigma_2)) < \epsilon$; 

\item[(c)] Letting $\Sym (n)$ act on $G_i ^n$ by permuting co-ordinates, 
$$\psi (\sigma_1) \tau((g_x)_x) \psi (\sigma_1)^{-1}
= \tau ((g_{\sigma_1 \cdot x})_x).$$
\end{itemize}

We shall call $\tau$ a \emph{base map} and $\psi$
an \emph{acting map}.
\end{Def}

\subsection{Sofic groups}

For the finite symmetric groups, we shall use both the notation 
$\Sym(A)$, to refer to the group of permutations 
of a specific finite set $A$, and the notation $\Sym(n)$, 
to refer to the group of permutations of the set $[n]=[1,n] \cap \mathbb{N}$. 

\begin{Def}\label{HammingDef}
Let $A$ be a finite set. The \emph{normalized Hamming metric}
on $\Sym(A)$ is given by:
$$d_A (\sigma, \tau) = 1 - \frac{1}{\lvert A \rvert} \lvert \lbrace a \in A : \sigma (a) = \tau(a) \rbrace \rvert$$
for $\sigma, \tau \in \Sym(A)$. 
In the case $A = [n]$, we may write $d_n$ for $d_A$. 
Let $\mathcal{G}_{\sof}$ denote the metric family of groups 
$\lbrace (\Sym(n),d_n) : n \in \mathbb{N} \rbrace$. 
A group is \emph{sofic} if it is $\mathcal{G}_{\sof}$-approximable. 
\end{Def}

\begin{Ex}
Every residually finite group and every amenable group is sofic. 
\end{Ex}

\begin{Lemma} \label{SoficUnifSepLemma}
Every sofic group  is $(\mathcal{G}_{\sof},1/5)$-approximable. 
\end{Lemma}

\begin{proof}
Let $\Gamma$ be a sofic group. Let $F \subseteq \Gamma$ be finite 
and let $\epsilon > 0$. We may assume $e \in F$ and $\epsilon < 1/10$. 
By Theorem 3.5 of \cite{Pestov}, for some $n \in \mathbb{N}$ there is a function 
$\varphi : \Gamma \rightarrow \Sym(n)$ with respect to $c=1/4$ and $\delta = \epsilon/2$. 
The conclusion is then immediate from Lemma \ref{BoringUnitalLemma}. 
\end{proof}

\begin{Lemma} \label{SoficCompatLem}
    The metric family $\mathcal{G}_{\sof}$ is product-compatible 
    and wreath-compatible. 
\end{Lemma}

\begin{proof}
    For finite sets $A$ and $B$, 
    define $\Delta_{A,B} : \Sym(A) \times \Sym(B) \rightarrow \Sym(A \times B)$ by: 
    $$\Delta_{A,B} (\sigma,\tau)[(a,b)] = (\sigma(a),\tau(b)).$$
    Then $(a,b) \in A \times B$ is a fixed point of 
    $\Delta_{A,B} (\sigma,\tau)$ iff $a$, $b$ are fixed points 
    of $\sigma$, $\tau$, respectively. 
    Thus for $\sigma_1 , \sigma_2 \in \Sym(A)$ and 
    $\tau_1 , \tau_2 \in \Sym(B)$, 
    with $d_A(\sigma_1 , \sigma_2) = x$ and $d_B(\tau_1 , \tau_2) = y$, we have:
    $$d_{A\times B} (\Delta_{A,B}(\sigma_1,\tau_1),\Delta_{A,B}(\sigma_2,\tau_2) ) = x+y-xy$$
    from which product-compatibility easily follows. 

    For wreath-compatibility, for $A$ a finite set 
    and $n \in \mathbb{N}$ define the homomorphisms 
    $\tau : \Sym(A)^n \rightarrow \Sym(A^n)$ 
    and $\psi : \Sym(n) \rightarrow \Sym(A^n)$ by: 
$$\tau ((\sigma_x)_{1\leq x\leq n}) [(a_y)_{1\leq y\leq n}]
= (\sigma_x(a_x))_{1\leq x\leq n}$$
    and: 
    $$\psi(\theta)[(a_x)_{1\leq x\leq n}] = (a_{\theta(x)})_{1\leq x\leq n}. $$
    Then $(a_x)_{1\leq x\leq n}$ is a fixed-point of 
    $\tau ((\sigma_x)_{1\leq x\leq n})$ iff for all 
    $1 \leq y \leq n$, $a_y$ is a fixed-point of $\sigma_y$, 
    and is a fixed-point of $\psi(\theta)$ iff 
    for all $1 \leq x \leq n$ $a_{\theta(x)} = a_x$, 
    whence wreath-compatibility easily follows. 
\end{proof}

\subsection{Linear sofic groups}

Let $\mathbb{K}$ be a field, and let $\GL_n(\mathbb{K})$ 
be the $n$-dimensional general linear group over $\mathbb{K}$. 

\begin{Def}
    The \emph{normalized rank-metric} $d_{\rk,n}$ on 
    $\GL_n(\mathbb{K})$ is given by 
    $$d_{\rk,n}(M,N) = \frac{1}{n}\rk(M-N).$$
    We write $\mathcal{G}_{\mathbb{K}}$ 
    for the metric family of groups $(\GL_n(\mathbb{K}),d_{\rk,n})_n$. 
\end{Def}

\begin{Rem}
        $d_{\rk,n}$ is clearly a bi-invariant metric on 
    $\GL_n(\mathbb{K})$. 
\end{Rem}

\begin{Def}
    The group $\Gamma$ is \emph{$\mathbb{K}$-linear sofic} 
    if it is $\mathcal{G}_{\mathbb{K}}$-approximable. 
\end{Def}

The class of linear sofic groups was introduced in \cite{ArzhPaunLinear}. 

\begin{Lemma} \label{LinSoficUnifSepLemma}
For every field $\mathbb{K}$ and every $c \in (0,1/8)$, 
every $\mathbb{K}$-linear sofic group is $(\mathcal{G}_{\mathbb{K}},c)$-approximable. 
\end{Lemma}

\begin{proof}
    This is immediate from Theorem 3.13 in \cite{Stolz}, using Lemma \ref{BoringUnitalLemma}.  
    (the case $\mathbb{K} = \mathbb{C}$ previously appeared as Proposition 5.13 of 
    \cite{ArzhPaunLinear}). 
\end{proof}

\begin{Lemma} \label{LinSoficCompatLemma}
    Let $\mathbb{K}$ be a field. Then the family $\mathcal{G}_{\mathbb{K}} = (\GL_n(\mathbb{K}),d_{\rk,n})_n$ is product-compatible and wreath-compatible.
\end{Lemma}

For the sake of the proof, it is convenient to introduce the 
following bi-invariant \emph{pseudo-metric} on $\GL_n(\mathbb{K})$: 
$$d_{\overline{\rk},n} (M,N) = \frac{1}{n} \min \lbrace \rk(M- \lambda N) : \lambda \in \mathbb{K}^{\ast} \rbrace$$
so that $d_{\overline{\rk},n} (M,N) \leq d_{\rk,n}(M,N)$, 
with equality iff $1 \in \mathbb{K}$ has maximal geometric multiplicity 
among the eigenvalues of $MN^{-1}$. 

Recall that taking the tensor-product of two finite-dimensional vector spaces 
$\mathbb{K}^m$ and $\mathbb{K}^k$ induces an injective homomorphism 
$(M,N) \mapsto M \otimes N$ from $\GL_m(\mathbb{K}) \times \GL_k (\mathbb{K})$ 
to $\GL_{mk} (\mathbb{K})$. 
The following observation from \cite{HayesSaleMetric} (Proposition 4.11) 
will be useful for the proof of product-compatibility.

\begin{Prop} \label{LinSoficProdBound}
    For all $M \in \GL_m(\mathbb{K})$ and $N \in \GL_k(\mathbb{K})$ we have: 
    $$d_{\overline{\rk},mk} (M\otimes N,I_{mk}) 
    \geq \max (d_{\overline{\rk},m}(M,I_m),d_{\overline{\rk},k}(N,I_k)).$$
\end{Prop}

\begin{proof}[Proof of Lemma \ref{LinSoficCompatLemma}]
We first prover wreath-compatibility. Let $m,k \in \mathbb{N}$, $n=km$
and identify $n \times n$ matrices $\mathbb{M}_n (\mathbb{K})$
over $\mathbb{K}$ with $k \times k$ block-matrices
$\mathbb{M}_k (\mathbb{M}_m (\mathbb{K}))$ over
$\mathbb{M}_m (\mathbb{K})$.
Let $\tau : \GL_m (\mathbb{K}) ^k \rightarrow \GL_n (\mathbb{K})$
be given by:
\begin{center}
$\tau((M_i)_i) = \left( \begin{array}{cccc}
M_1 & 0 & \cdots & 0 \\
0 & M_2 & \cdots & 0 \\
\vdots & \vdots & \ddots & \vdots \\
0 & 0 & \cdots & M_k
\end{array} \right)$
\end{center}
and let $\psi : \Sym(k) \rightarrow \GL_n (\mathbb{K})$
be given by $\psi(\sigma) = (A_{i,j})_{i,j=1} ^k$,
where $A_{i,j} \in \mathbb{M}_m (\mathbb{K})$ is given by:
\begin{center}
$A_{i,j} = \begin{cases}
I_m & i=\sigma(j) \\
0 & \text{otherwise}
\end{cases} $. 
\end{center}
Then $\tau$ and $\psi$ are homomorphisms and we have:
$$\frac{1}{n} \rk (\tau((M_i)_i)-I_n) = \frac{1}{k} \sum_{j=1} ^k \frac{1}{m} \rk  (M_j - I_m)$$
while:
$$\frac{1}{n} \rk (\psi(\sigma)-I_n) = n - m \cyc (\sigma) \leq m \lvert \supp (\sigma) \rvert$$
where $\cyc(\sigma)$ is the number of cycles in a disjoint
cycle decomposition of $\sigma$ (including $1$-cycles),
so that $\cyc(\sigma) \geq k - \lvert \supp(\sigma) \rvert$.
It follows immediately that $\tau$ and $\psi$ satisfy the defining properties 
of a base and acting map, respectively. 

We now prove product-compatibility. For $m \in \mathbb{N}$ and $A \in \mathbb{M}_m (\mathbb{K})$
define $\hat{A} \in \mathbb{M}_{2m} (\mathbb{K})$ by:
\begin{center}
$\hat{A} = \left( \begin{array}{cc}
A & 0 \\
0 & I_m
\end{array} \right)$
\end{center}
so that: $$d_{\overline{\rk},2m} (\hat{A},I_{2m})
= d_{\rk,2m} (\hat{A},I_{2m})
= \frac{1}{2} d_{\rk,m} (A,I_m). $$
For $k,m\in\mathbb{N}$, define $\Delta_{m,k} : \GL_m (\mathbb{K}) \times \GL_k (\mathbb{K})
\rightarrow \GL_{4km} (\mathbb{K})$
by $$\Delta_{m,k}(A,B) = \hat{A} \otimes \hat{B}.$$
Then $\Delta_{m,k}$ is a homomorphism.
Suppose that $A$ and $B$ have $1$-eigenspaces of dimensions
$m_1$ and $k_1$, respectively.
Then the $1$-eigenspace of $\Delta_{m,k}(A,B)$ has dimension
at least $(m+m_1)(k+k_1)$, so if
$d_{\rk,m} (A,I_m) = \delta$ and $d_{\rk,k} (B,I_k) = \epsilon$,
we have $$d_{\rk,4mk} ( \Delta_{m,k}(A,B),I_{4mk} ) \leq \frac{2\delta+2\epsilon-\delta\epsilon}{4} \leq \frac{\delta + \epsilon}{2}.$$
On the other hand, $\Delta_{m,k}(A,B)$ always has $1$ eigenspace
of dimension at least $mk$.
If $1$ has maximal geometric multiplicity
among the eigenvalue of $\Delta_{m,k}(A,B)$,
then using Proposition \ref{LinSoficProdBound} we have:
\begin{align*}
d_{\rk,4mk} (\Delta_{m,k}(A,B),I_{4mk})
& = d_{\overline{\rk},4mk} (\hat{A} \otimes \hat{B},I_{4mk}) \\
& \geq \max \lbrace d_{\overline{\rk},2m} (\hat{A},I_{2m}),d_{\overline{\rk},2k} (\hat{B},I_{2k}) \rbrace \\
& = \frac{1}{2} \max \lbrace d_{\rk,m} (A,I_m),d_{\rk,k} (B,I_k) \rbrace
\end{align*}
If $1$ does \emph{not} have maximal geometric multiplicity
among the eigenvalue of $\Delta_{m,k}(A,B)$,
then $d_{\rk,4mk} (\Delta_{m,k}(A,B),I_{4mk}) \geq 1/4$. Thus in any case we have: 
$$d_{\rk,4mk} (\Delta_{m,k}(A,B),I_{4mk}) \geq \frac{1}{4} \max \lbrace d_{\rk,m} (A,I_m),d_{\rk,k} (B,I_k).$$
Hence $\Delta_{m,k}$ satisfies the defining properties of a product map. 
\end{proof}

\subsection{Hyperlinear groups}

Let $U(n) \subseteq \mathbb{M}_n (\mathbb{C})$ denote the 
group of $n \times n$ unitary matrices. 

 \begin{Def}
     The \emph{Hilbert-Schmidt} metric $d_{\HS_n}$ on $U(n)$ 
     is given by: 
     $$d_{\HS_n} (A,B) = \sqrt{\frac{1}{n} \sum_{i,j=1} ^n \lvert A_{i,j} - B_{i,j} \rvert^2} = \sqrt{\frac{1}{n} \Tr((A-B)(A-B)^{\ast})}. $$
 \end{Def}

 It is easy to see that $d_{\HS_n}$ is a bi-invariant metric 
 on $U(n)$. 

 \begin{Def}
     We denote by $\mathcal{G}_{\hyp}$ 
     the metric family of groups 
     $\lbrace (U(n),d_{\HS_n}) : n \in \mathbb{N} \rbrace$. 
     A group is described as \emph{hyperlinear} 
     if it is $\mathcal{G}_{\hyp}$-approximable. 
 \end{Def}

 \begin{Rem} \label{HammingHSRmrk}
     For $\sigma , \tau \in \Sym(n)$, the associated permutation matrices $P(\sigma) , P(\tau) \in \mathbb{M}_n (\mathbb{C})$ 
     are unitary, and satisfy: 
     $$d_n(\sigma , \tau) = \frac{1}{2} d_{\HS_n} (P(\sigma) , P(\tau))^2.$$
     Thus every sofic group is hyperlinear. 
 \end{Rem}

\begin{Lemma} \label{HypUnifSepLemma}
Every hyperlinear group is $(\mathcal{G}_{\hyp},1/5)$-approximable. 
\end{Lemma}

\begin{proof}
    This is identical to the proof of Lemma \ref{SoficUnifSepLemma}, 
    but using Theorem 3.6 from \cite{Pestov} instead of Theorem 3.5. 
\end{proof}

\begin{Lemma} \label{HypCompatLem}
    The metric family $\mathcal{G}_{\hyp}$ is product-compatible 
    and wreath-compatible. 
\end{Lemma}

Similarly to the case of linear sofic groups, it is convenient to introduce the 
following bi-invariant {pseudo-metric} on $U(n)$: 
$$d_{\overline{\HS},n} (M,N) = \min \lbrace d_{{\HS},n}(M- \lambda N) \mid \lambda \in U(1) \rbrace$$
so that $d_{\overline{\HS},n} (M,N) \leq d_{\HS,n}(M,N)$.

Similar to the linear sofic case, we then have

\begin{Prop} \label{ProjHSProdBound}
    For all $M \in U(m)$ and $N \in U(k)$ we have: 
    $$d_{\overline{\HS},mk} (M\otimes N,I_{mk}) 
    \geq \max (d_{\overline{\HS},m}(M,I_m),d_{\overline{\HS},k}(N,I_k)).$$
\end{Prop}

\begin{proof}[{Proof of Lemma \ref{HypCompatLem}}]
    For wreath-compatibility, 
    let $k,m \in \mathbb{N}$m set $n=km$ and 
    define the homomorphisms 
    $\tau : U(m)^k \rightarrow \GL_n (\mathbb{C})$ 
    and $\psi : \Sym(k) \rightarrow \GL_n (\mathbb{C})$
    as in the proof of Lemma \ref{LinSoficCompatLemma}. 
    Then $\im (\tau) , \im (\psi) \leq U(n)$. 
    Moreover, for $A_i , B_i \in U(m)$ ($1 \leq i \leq k$), 
    we have: 
    $$n d_{\HS_n} (A,B)^2 = m \sum_{i=1} ^k d_{\HS_m} (A_i,B_i)^2$$
    while for $\sigma \in \Sym(k)$, 
    $\psi(\sigma) \in U(n)$ is the permutation matrix 
    $P(\sigma')$ of a permutation $\sigma' \in \Sym(n)$ 
    satisfying $d_n(\sigma' \id_n) = d_k(\sigma \id_k)$. 
    Thus $\tau$ and $\psi$ are seen to satisfy the defining 
    properties of a base and acting map 
    (using the equation in Remark \ref{HammingHSRmrk} 
    for $\psi$).     

    For product-compatibility, 
    by analogy with relevant part of the proof of Lemma \ref{LinSoficCompatLemma} , 
    we again define for $A\in U(m)$ the matrix
\begin{center}
$\hat{A} = \left( \begin{array}{cc}
A & 0 \\
0 & I_m
\end{array} \right) \in U(2m)$
\end{center}
and notice that
\[
d_{\overline{\HS},2m}(\hat A,I_{2m})\geq \frac{1}{4} d_{\HS,m}(A,I_m).
\]
For $k,m\in\mathbb{N}$, define $\Delta_{m,k} : \GL_m (\mathbb{K}) \times \GL_k (\mathbb{K})
\rightarrow \GL_{4km} (\mathbb{K})$
by $$\Delta_{m,k}(A,B) = \hat{A} \otimes \hat{B}.$$
Then $\Delta_{m,k}$ is a homomorphism and if
$d_{\HS,m} (A,I_m) = \delta$ and $d_{\HS,k} (B,I_k) = \epsilon$,
we have $$d_{\HS,4mk} ( \Delta_{m,k}(A,B),I_{4mk} ) \leq {\delta + \epsilon}.$$

On the other hand, $\Delta_{m,k}(A,B)$ always has $1$ eigenspace
of dimension at least $mk$. If $1$ has maximal geometric multiplicity
among the eigenvalue of $\Delta_{m,k}(A,B)$,
then using Proposition \ref{ProjHSProdBound} we have:
\begin{align*}
d_{\HS,4mk} (\Delta_{m,k}(A,B),I_{4mk})
& = d_{\overline{\HS},4mk} (\hat{A} \otimes \hat{B},I_{4mk}) \\
& \geq \max \lbrace d_{\overline{\HS},2m} (\hat{A},I_{2m}),d_{\overline{\HS},2k} (\hat{B},I_{2k}) \rbrace \\
& \geq \frac{1}{4} \max \lbrace d_{\HS,m} (A,I_m),d_{\HS,k} (B,I_k) \rbrace
\end{align*}
If $1$ does \emph{not} have maximal geometric multiplicity
among the eigenvalue of $\Delta_{m,k}(A,B)$,
then $d_{\HS,4mk} (\Delta_{m,k}(A,B),I_{4mk}) \geq 1/4$. Thus in any case we have: 
$$d_{\HS,4mk} (\Delta_{m,k}(A,B),I_{4mk}) \geq \frac{1}{4} \max \lbrace d_{\HS,m} (A,I_m),d_{\HS,k} (B,I_k).$$
Hence $\Delta_{m,k}$ satisfies the defining properties of a product map. 

\end{proof}

\subsection{Weakly sofic groups}

The class of weakly sofic groups was introduced in \cite{GlebRive}. 
Let $\mathcal{G}_{\weak}$ denote the metric family 
of all finite groups equipped with bi-invariant metrics 
of diameter one. 

\begin{Def}
    A group is described as \emph{weakly sofic} if it is 
    $\mathcal{G}_{\weak}$-approximable. 
\end{Def}

\begin{Rem}
    It is immediate from the definitions that 
    sofic groups are weakly sofic. 
    For any field $\mathbb{K}$, 
    every $\mathbb{K}$-linear sofic group 
    is weakly sofic (Theorem 1.3 of \cite{ArzhPaunLinear}). 
\end{Rem}

\begin{Lemma} \label{WSoficUnifSepLemma}
    For every $c \in (0,1)$, every weakly sofic group is $(\mathcal{G}_{\weak},c)$-approximable. 
\end{Lemma}

\begin{proof}
Let $d$ be a bi-invariant metric of diameter one on a finite group $G$. 
    Let $f : [0,1] \rightarrow [0,1]$ be given by $f(x)=2x-x^2$. 
    Then $f$ is an increasing subadditive bijection, 
    so that $f \circ d$ is also a bi-invariant metric of diameter one. 
Moreover, for all $x \in [0,1]$ we have $x leq f(x) \leq 2x$ and 
if $x \leq a \in (0,1)$ then $f (x) \geq (2-c)x$. 

Let $\Gamma$ be a weakly sofic group, and suppose that 
$c' > 0$ is such that $\Gamma$ is $(\mathcal{G}_{\weak},c')$-approximable. 
Let $n \in \mathbb{N}$ be such that $(2-c)^n c' \geq c$. 
Let $F \subseteq \Gamma$ be a finite set, let $\epsilon > 0$ 
and let $\varphi : \Gamma \rightarrow G$ be an $(F,\epsilon/2^n,c')$-$\mathcal{G}_{\weak}$-representation of $\Gamma$, for some finite group $G$ equipped with 
a bi-invariant metric of diameter one. 
Then equipping $G$ with the metric $f^n \circ d$ instead of $d$, 
$\varphi$ becomes an $(F,\epsilon,c')$-$\mathcal{G}_{\weak}$-representation of $\Gamma$. 
\end{proof}

\begin{Lemma} \label{WSoficCompatLem}
    The metric family $\mathcal{G}_{\weak}$ 
    is product-compatible and wreath-compatible. 
\end{Lemma}

\begin{proof}
For product compatibility, 
if $(G_1,d_1),(G_2,d_2) \in \mathcal{G}_{\weak}$ 
set $G = G_1 \times G_2$ and let 
$d : G \times G \rightarrow [0,1]$ be given by: 
$$d((g_1,g_2),(g_1',g_2')) = \max \lbrace d_1 (g_1,g_1 ') , 
d_2 (g_2,g_2') \rbrace$$
so that $d$ is a bi-invariant metric on $G$ of diameter one, 
and the identity map $G_1 \times G_2 \rightarrow G$ 
witnesses product compatibility. 

For wreath-compatibility, given $(G,d) \in \mathcal{G}_{\weak}$ 
and $n \in \mathbb{N}$, 
let $G'$ be the permutational wreath product 
$G \wr \Sym(n) = G^n \rtimes \Sym(n)$ and 
define $d' : G' \times G' \rightarrow [0,1]$ by: 
$$d'\Big( \big((g_i)_i,\sigma\big) ,  \big((g' _i)_i,\sigma'\big)  \Big) 
= d_n (\sigma,\sigma')+\frac{1}{n} \sum_{\substack{1 \leq i \leq n \\ \sigma(i)=\sigma'(i)}} d_i (g_i,g' _i)$$
Then $d'$ is easily seen to be a metric on $G'$ of diameter one, 
and by Proposition 2.9 of \cite{HayesSaleMetric} it is bi-invariant. 
Then $\tau : G^n \rightarrow G'$ and $\psi : \Sym(n) \rightarrow G'$, 
given by $\tau (((g_i)_i))=\big((g_i)_i,\id_n \big)$ and 
$\psi (\sigma) = ((e)_i,\sigma)$ satisfy the conditions of 
a base map and acting map, respectively. 
\end{proof}

\subsection{LEF groups}

The class of LEF groups (that is, groups \emph{locally embeddable into the class of finite groups}) 
was introduced by Vershik and Gordon \cite{VersGord}. 

\begin{Def} \label{LEFGrpDef}
    For $\Gamma$ and $\Delta$ groups and $A \subseteq \Gamma$, 
    a function $\phi: A \rightarrow \Delta$ is a 
    \emph{partial homomorphism} if, for all $g,h \in A$, 
    if $gh \in A$ then $\phi(gh)=\phi(g)\phi(h)$. 
    For $\mathcal{C}$ a class of groups, closed under taking 
    subgroups, the group $\Gamma$ is 
    \emph{locally embeddable into $\mathcal{C}$} if, 
    for every finite subset $A \subseteq \Gamma$, there 
    is a group $\Delta$ lying in $\mathcal{C}$ 
    and an injective partial homomorphism 
    $\phi : A \hookrightarrow \Delta$. 
    In the case for which $\mathcal{C}$ is the class of all 
    finite groups, a group which is locally embeddable into $\mathcal{C}$ 
    shall also be called an \emph{LEF group}. 
\end{Def}

It is also possible to view the class of LEF groups 
as those approximable in a metric class, in the sense 
of Subsection \ref{GenMetricApproxSubsect}. 
Namely, we take $\mathcal{G}_{\LEF}$ the class of all finite groups $G$ 
equipped with the discrete metric 
$d : G \times G \rightarrow \lbrace 0,1 \rbrace$ 
given by $d(g,h) = \delta_{g,h}$. 
The following is then an easy consequence of definitions. 

\begin{Lemma} \label{LEFMetricLem}
    The group $\Gamma$ is LEF iff it is $\mathcal{G}_{\LEF}$-approximable. 
\end{Lemma}

\begin{Def}
    Suppose the class $\mathcal{C}$ is closed under taking subgroups 
    and finite direct products. The group $\Gamma$ is 
    \emph{residually in $\mathcal{C}$} if, for every finite subset $A \subseteq \Gamma$, 
    there is a function $\phi$ as in 
    Definition \ref{LEFGrpDef}, which is the restriction of a homomorphism 
    $\Gamma \rightarrow \Delta$. 
\end{Def}

Thus, for example, every residually finite group is LEF. 
Vershik and Gordon note that, within the class of finitely presented groups, 
the two concepts coincide \cite{VersGord}. 

\begin{Prop} \label{fpLEFRFProp}
    Every finitely presented LEF group is residually finite. 
\end{Prop}

\begin{Rem} \label{LELEFRem}
If $\mathcal{G}$ is a metric family of groups and $\mathcal{C}$ is the class of 
$\mathcal{G}$-approximable groups, then every group which is locally 
embeddable into $\mathcal{C}$ lies in $\mathcal{C}$. 
In particular, every group which is locally embeddable into the 
class of LEF groups is LEF. 
\end{Rem}

Finally, the following closure properties of the class of LEF groups are useful. 

\begin{Thm} \label{VersGordWPLEF}
The restricted regular wreath product of two LEF groups is LEF. 
\end{Thm}

\begin{Prop} \label{VLEFisLEFProp}
    Suppose $\Gamma$ has a finite-index subgroup $\Delta$ which is LEF. 
    Then $\Gamma$ is LEF. 
\end{Prop}

\begin{proof}
    Passing to a further finite-index subgroup we may assume $\Delta \vartriangleleft \Gamma$ 
    (since the class of LEF groups is subgroup-closed). 
    By the Kaluzhnin--Krasner Theorem, $\Gamma$ embeds as a subgroup of 
    $\Delta \wr (\Gamma/\Delta)$ 
    (since $\Gamma/\Delta$ is finite, the restricted and unrestricted wreath products coincide). 
    The latter is a LEF group, since $\Delta$ and $\Gamma/\Delta$ are, 
    by Theorem \ref{VersGordWPLEF}. 
\end{proof}

\section{Approximations of actions}

\subsection{LEF actions on sets and graphs}

The following notion appears in \cite{BradLEFExtn}. 

\begin{Def} \label{LEFActionDef}
A group action $\Gamma \curvearrowright_{\alpha} X$ 
on a set $X$ 
is a \emph{LEF action} if,
for all finite subsets $F \subseteq \Gamma$ and $Z \subseteq X$,
there exists a finite group $Q$, a finite set $Y$,
an action $Q \curvearrowright_{\beta} Y$ and functions
$\rho : \Gamma \rightarrow Q$, $\pi : Z \rightarrow Y$,
with $\pi$ injective, such that:
\begin{itemize}
\item[(i)] For all $g,h \in F$, if $gh \in F$ then $\rho(gh)=\rho(g)\rho(h)$;

\item[(ii)] If $g \in F$ and $x \in Z$ are such that $\alpha(g)\cdot x \in Z$,
then $\pi(\alpha(g)\cdot x) = \beta(\rho(g)) \cdot \pi(x)$.
\end{itemize}
If $X=(V,E)$ has the structure of a graph and $\Gamma \curvearrowright_{\alpha} X$ is an action by graph automorphisms, 
we describe $\alpha$ as a \emph{LEF action} if, 
for all finite subsets $F \subseteq \Gamma$ and 
finite induced subgraphs $Z \subseteq X$,
there exists a finite group $Q$, a finite graph $Y$,
an action $Q \curvearrowright_{\beta} Y$ by graph automorphisms 
and functions $\rho,\pi$ as above, 
with $\pi$ an inclusion of an induced subgraph. 
\end{Def}

In the sequel, whenever LEF actions arise it shall be clear from 
context whether we intend a LEF action in the 
set-theoretic or graph-theoretic category. 

\begin{Rem} \label{LEFActionBasicsRem}
    The following properties of LEF actions on sets (respectively graphs) 
    are easy consequences of Definition \ref{LEFActionDef}. 
    \begin{itemize}
        \item[(i)] If $\phi : \Delta \rightarrow \Gamma$ is a homomorphism, 
        and $\Gamma \curvearrowright_{\alpha} X$ is a LEF action on a set (graph), 
        then the induced action $\Delta \curvearrowright_{(\alpha \circ \phi)} X$ 
        is a LEF action. In particular, the restriction of a LEF action to a 
        subgroup is LEF; 

        \item[(ii)] If $Y \subseteq X$ is an invariant subset (induced subgraph) 
        under the LEF action $\Gamma \curvearrowright_{\alpha} X$, 
        then the restriction of $\alpha$ to $Y$ is a LEF action; 

        \item[(iii)] If $\Gamma = \bigcup_{i \in I} \Gamma_i$ is an ascending union of subgroups 
        $\Gamma_i$, and $\Gamma \curvearrowright_{\alpha} X$ is an action whose 
        restriction to each $\Gamma_i$ is a LEF action, then $\alpha$ is a LEF action 
        of $\Gamma$; 

        \item[(iv)] If $X$ is a set, and the restriction of the action $\Gamma \curvearrowright_{\alpha} X$ to each orbit of $\alpha$ is LEF, then $\alpha$ is LEF; 

        \item[(v)] If $X=(V,E)$ is a graph and $\Gamma \curvearrowright_{\alpha} X$ 
        is a LEF action, then the action of $\Gamma$ on the set $V$ and the 
        graph-complement $X^c$ are LEF actions. 
        Conversely, if $V$ is a set admitting a LEF action of $\Gamma$, 
        then the actions of $\Gamma$ on the both complete graph $K(V)$ on $V$ 
        and the empty graph $K(V)^c$ on $V$ are LEF actions. 
        
    \end{itemize}
\end{Rem}

It turns out that, for transitive actions, a concept equivalent to that of a LEF action 
was independently introduced by Cornulier \cite{CornulierWreath}. 
Recall that for a group $\Gamma$, $\Sub(\Gamma)$ denotes the 
space of subgroups of $\Gamma$, viewed as a (closed) subspace of 
$\lbrace 0,1 \rbrace ^{\Gamma}$ (equipped with the product topology). 
Thus, for $L , L_n \leq \Gamma$ ($n \in \mathbb{N}$), 
the sequence $(L_n)$ converges to $L$ in $\Sub(\Gamma)$ 
iff, for all finite subsets $F \subseteq \Gamma$, 
$L_n \cap F$ eventually equals $L \cap F$. 

\begin{Def}
A subgroup $\Delta$ of the groups $\Gamma$ is almost 
of finite index (AFI) if, for every finitely presented group 
$\tilde{\Gamma}$ and every homomorphism 
$\rho : \tilde{\Gamma} \rightarrow \Gamma$, 
there is a sequence of finite-index subgroups 
of $\tilde{\Gamma}$ converging to $\rho^{-1} (\Gamma)$ 
in $\Sub(\tilde{\Gamma})$. 
\end{Def}

\begin{Rem} \label{AFIRmrk}
It is easy to see (comment immediately following Definition 1.1 of \cite{CornulierProper2012}) 
that the set $\Sub_{\AFI}(\Gamma) \subseteq \Sub(\Gamma)$ of AFI subgroups 
of $\Gamma$ forms a closed subspace of $\Sub(\Gamma)$. 
\end{Rem}

\begin{Prop} \label{AFILEFActProp}
Let $\Delta \leq \Gamma$. 
Then $\Delta$ is AFI in $\Gamma$ iff the quasiregular action 
$\Gamma \curvearrowright \Gamma/\Delta$ is a LEF action. 
\end{Prop}

\begin{proof}
First suppose $\Gamma \curvearrowright_{\alpha} X$ 
is a transitive LEF action and $x_0 \in X$ satisfying 
$\Stab_{\Gamma} (x_0) = \Delta$. 
Let $\tilde{\Gamma}$ be a finitely presented group; 
$\rho : \tilde{\Gamma} \rightarrow \Gamma$ be a homomorphism, 
and let $B \subseteq \tilde{\Gamma}$ be a finite subset. 
Fix a finite generating set $S$ for $\tilde{\Gamma}$, 
and let $C>0$ be such that $\tilde{\Gamma}$ 
admits a presentation with generating set $S$ and all relations 
of word-length at most $C$ in $S$. 
Let $F \subseteq \Gamma$ be a finite set containing both $\rho(B)$ 
and the image under $\rho$ of all words in $S$ of length at most $C$. 
In particular $e \in F$. 
Let $Z \subseteq X$ be a finite subset containing 
$\lbrace \alpha(g)\cdot x_0 : g \in F \rbrace$. 
Associated to $F$ and $Z$ we have $\psi : F \rightarrow Q$; 
$\pi : Z \hookrightarrow Y$ and $Q \curvearrowright_{\beta} Y$ 
as in the definition of LEF action. 
By choice of $F$, $\psi$ induces a well-defined homomorphism 
$\tilde{\psi} : \tilde{\Gamma} \rightarrow Q$, 
hence an induced action $\tilde{\Gamma} \curvearrowright_{\tilde{\beta}} Y$. 
Let $y_0 = \pi(x_0)$ and set $L = \Stab_{\tilde{\Gamma}} (y_0)$, 
a finite-index subgroup of $\tilde{\Gamma}$. 
We claim that $L \cap B = \rho^{-1} (\Delta) \cap B$, 
which completes the proof. 
This is so, because for $g \in B$, $\rho(g) \in F$ and 
$x_0 , \rho (g) \cdot x_0 \in Z$ with 
$\pi(\rho (g) \cdot x_0) = \psi (\rho(g))\cdot y_0$, so that 
$g \in \rho^{-1} (\Delta)$ iff $\rho (g) \cdot x_0 = x_0$ 
iff $\psi(\rho(g))\cdot y_0 = y_0$. 

Conversely suppose that $\Delta \leq \Gamma$ is an AFI subgroup; 
let $X = \Gamma/\Delta$; let $\Gamma \curvearrowright_{\alpha} X$ 
be the quasiregular action, and let $F \subseteq \Gamma$, 
$Z \subseteq X$ be finite. 
Expanding $F$, we can assume that $F$ contains a 
set $\lbrace g_1 , \ldots , g_{\lvert Z \rvert} \rbrace$ 
of coset representatives for the elements of $Z$, 
and that $e = F = F^{-1}$. 
Let $\Gamma' = \langle F \rangle \leq \Gamma$. 
Pick a presentation $\Gamma' \cong \langle S | R \rangle$ for $\Gamma'$, 
with $S \subseteq \Gamma'$ finite. 
Let $l_0 > 0$ be such that $F \subseteq B_S(l_0)$. 
Let $F(S)$ be the free group on basis $S$, 
and for $l > l_0$ let $R_l \subseteq F(S)$ be given by $R_l = \ll R \gg^{F(S)} \cap B_S (l)$; 
let $\Gamma_l ' = \langle S | R_l \rangle$, and let 
$\rho_l : \Gamma_l ' \twoheadrightarrow \Gamma'$ be the natural epimorphism. 
Then there exists $l_1 > 3 l_0$ such that $\rho_{l_1} |_{B_S (3l_0)}$ is injective. 
Set $\tilde{\Gamma} = \Gamma_{l_1}$ (a finitely presented group), and take $\rho = \rho_{l_1}$. 
There is a unique injection $\tilde{\psi} : B_S (3l_0) \rightarrow \tilde{\Gamma}$ 
(defined on $B_S (3l_0) \subseteq \Gamma'$) 
such that $\rho \circ \tilde{\psi} = \id_{B_S (3l_0)}$. 
It is clear that $\tilde{\psi}$ is a partial homomorphism. 

Applying the AFI hypothesis, there exists $L \leq_f \tilde{\Gamma}$ 
such that $B_S (l_1) \cap L = B_S (l_1) \cap \rho^{-1} (\Delta)$ in $\tilde{\Gamma}$. 
Let $N \vartriangleleft \tilde{\Gamma}$ be the normal core of $L$ in $\tilde{\Gamma}$, 
and set $Q = \tilde{\Gamma} / N$ (a finite group) and $Y = \tilde{\Gamma}/L$, 
so that we have a natural action $Q \curvearrowright_{\beta} Y$ by left-multiplication 
of cosets, and $\tilde{\psi}|_F$ descends to $\psi : F \rightarrow Q$. 
We may define $\pi : Z \rightarrow Y$ by $\pi (g_i \Delta) = \tilde{\psi} (g_i) L$, 
for $1 \leq i \leq \lvert Z \rvert$. If $\tilde{\psi} (g_i) L = \tilde{\psi} (g_j) L$, 
then applying $\rho$, we have $g_j ^{-1} g_i \in \Delta$, so $\pi$ is injective. 
Moreover, if $g \in F$ and $x = g_i \Delta \in Z$, 
with $gx \in Z$, so that for some $1 \leq j \leq \lvert Z \rvert$ 
we have $g_j ^{-1} g g_i \in \Delta \cap B_S (3l_0)$, so that:
$$\tilde{\psi}(g_j)^{-1} \tilde{\psi} (g) \tilde{\psi} (g_i) 
= \tilde{\psi} (g_j ^{-1} g g_i) \in L$$
hence: 
$$\pi(gx) = \pi (g_j \Delta) = \tilde{\psi}(g_j) L = \tilde{\psi}(g)\tilde{\psi}(g_i)L = \beta(\psi(g)) [\pi(x)]$$
so the conditions of Definition \ref{LEFActionDef} are satisfied. 
\end{proof}




\begin{Prop} \label{LEFGrpFaithLEFActProp}
For a group $\Gamma$, the following are equivalent. 
\begin{itemize}
    \item[(i)] $\Gamma$ is a LEF group; 

    \item[(ii)] The action of $\Gamma$ on itself (by left multiplication) is a LEF action; 

    \item[(iii)] $\Gamma$ admits a faithful LEF action on a set. 
\end{itemize}
\end{Prop}

\begin{proof}
If $\Gamma$ is a LEF group then the action of $\Gamma$ 
on itself by left-multiplication satisfies the conditions 
of Definition \ref{LEFActionDef}: 
given $F,Z \subseteq \Gamma$ finite subsets, 
let $Q$ be a finite group such that there exists 
an injective partial homomorphism $\phi : F \cup Z \rightarrow Q$. 
Then $\rho : \Gamma \rightarrow Q$, 
obtained by arbitrarily extending $\phi|_F$ to $\Gamma$; 
$\pi = \phi|_Z : Z \hookrightarrow Y = Q$, 
and $Q$ acting on itself by left-multiplication, 
witness the conditions of Definition \ref{LEFActionDef}. 

The action of $\Gamma$ on itself (by left multiplication) is a faithful action, 
so (ii) implies (iii). 

Suppose $\Gamma \curvearrowright_{\alpha} X$ 
is a faithful LEF action on a set. 
Let $F \subseteq \Gamma$ be a finite subset. 
For each $e \neq g \in F$, let $x_g \in X$ with 
$\alpha(g)[x_g] \neq x_g$ and set 
$Z = \lbrace x_g , \alpha(g)[x_g] : g \in F \setminus \lbrace e \rbrace \rbrace$. 
Let $Q$, $Y$, $Q \curvearrowright_{\beta} Y$, 
$\rho : \Gamma \rightarrow Q$ and $\pi : Z \hookrightarrow Y$ 
be as in Definition \ref{LEFActionDef}. 
Then for $e \neq g \in F$, $x_g , \alpha(g)[x_g] \in Z$, so:  
$$\pi(x_g) \neq \pi(\alpha(g)[x_g]) = \beta(\rho(g))[\pi(x_g)].$$
Thus $\rho |_F$ is an injective partial homomorphism 
from $F$ to a finite group. 
Since this holds for every finite subset $F \subseteq \Gamma$, 
$\Gamma$ is a LEF group. 
\end{proof}

We recall that the group $\Gamma$ is \emph{locally extended residually finite (LERF)} 
is every finitely generated subgroup of $\Gamma$ is the intersection of finite-index 
subgroups of $\Gamma$. The next result is a slight sharpening of the 
main result from \cite{Gao}. 
 
\begin{Thm} \label{LERFLEFActThm}
    If $\Gamma$ is LERF, then all actions of $\Gamma$ on sets are LEF actions. 
\end{Thm}

\begin{proof}
    By Remark \ref{LEFActionBasicsRem} (iv), it is sufficient to consider 
    transitive actions. By Proposition \ref{AFILEFActProp}, 
it therefore suffices to prove that if $\Gamma$ is LERF, 
then every subgroup of $\Gamma$ is AFI. 
    This is indeed the case, as is explained in Example 2.6 of \cite{CornulierWreath}: 
    finitely generated subgroups of $\Gamma$ are dense in $\Sub(\Gamma)$, 
    then since $\Gamma$ is LERF, every finitely generated subgroup is a limit 
    in $\Sub(\Gamma)$ of finite index subgroups, and clearly finite-index subgroups 
    are AFI. 
\end{proof}

\begin{Ex}
    The following groups are known to be LERF, hence by Theorem \ref{LERFLEFActThm}, 
    all of their actions on sets are LEF actions. 
\begin{itemize}
    \item[(i)] All polycyclic groups (a 1948 result of Malcev, see for instance \cite{Alpe}); 

    \item[(ii)] All limit groups (including free and surface groups) \cite{Wilton}; 

    \item[(iii)] Fundamental groups of compact hyperbolic $3$-manifolds; 
    Sol $3$-manifolds, and Seifert fibre spaces (see references in Definition 2 of \cite{Gao}); 

    \item[(iv)] The regular restricted wreath product of two finitely generated abelian groups 
    is LERF, as are free metabelian groups \cite{AlpeWreath}. 
    
\end{itemize}
\end{Ex}

\subsection{LEF actions on groups}

In this Subsection we introduce an adaptation of Definition \ref{LEFActionDef} 
to the setting of an action of one group on another by automorphisms. 

\begin{Def} \label{LEFActGrpDef}
    Let $\Gamma$ and $\Delta$ be groups; let 
$\Gamma \curvearrowright_{\beta} \Delta$ 
be an action by automorphisms, and let $\mathcal{C}$ be a class of groups. 
The action $\beta$ is \emph{$\mathcal{C}$-LEF} if, 
for all finite subsets $F \subseteq \Gamma$ 
and $E \subseteq \Delta$, 
there is a finite group $Q$; a group $K \in \mathcal{C}$, 
an action $Q \curvearrowright_{\gamma} K$ by automorphisms, 
and partial homomorphisms 
$\rho: F \rightarrow Q$, $\pi: E \hookrightarrow K$, with $\pi$ injective, 
such that for all $g \in F$ and $h \in E$, 
if $\beta(g)[h] \in E$ then 
$\pi(\beta(g)[h]) = \gamma(\rho(g))[\pi(h)]$. 
\end{Def}

Note that, if there is a group admitting a $\mathcal{C}$-LEF action by automorphisms 
on the group $\Delta$ then $\Delta$ is locally embeddable into $\mathcal{C}$. 

\begin{Rem}
In the special case for which $\mathcal{C}$ is the class of finite groups 
and the partial homomorphism $\rho$ is assumed to be injective, 
we recover the notion of an \emph{equivariantly approximable} action, 
introduced in \cite{VersGord}. 
\end{Rem}

\begin{Ex}
    Let $\Delta$ be a finitely generated residually finite group. 
    A famous observation of Mal'cev yields that $\Aut(\Delta)$ is also 
    residually finite. By a similar argument, we see that any 
    action $\Gamma \curvearrowright_{\beta} \Delta$ of a group $\Gamma$ 
    on $\Delta$ by automorphisms is a LEF $\mathcal{F}$-action, 
    for $\mathcal{F}$ the class of finite groups. 
    For let $E \subseteq \Delta$ and $F \subseteq \Gamma$ be finite subsets. 
    Let $N \vartriangleleft_f \Delta$ be such that all elements of $E$ 
    lie in distinct cosets of $N$ in $\Delta$. 
    Since $\Delta$ is finitely generated, 
    $N$ contains a finite-index characteristic subgroup $M$ of $\Delta$. 
    The action of $\Gamma$ on $\Delta$ induces 
    an action of $\Gamma$ on $K = \Delta / M$; 
    let the kernel of this action be $L \vartriangleleft_f \Gamma$ and set $Q = \Gamma/L$. 
    The action $\Gamma \curvearrowright_{\beta} \Delta$ descends to an action 
    $Q \curvearrowright_{\gamma} K$, and letting $\rho : \Gamma \rightarrow Q$, 
    $\pi : \Delta \rightarrow K$ be the quotient homomorphisms, 
    the restrictions of $\rho$ and $\pi$ to $F$ and $E$ satisfy the conditions of 
    Definition \ref{LEFActGrpDef}. 
\end{Ex}

Vershik and Gordon \cite{VersGord} note that if $\Gamma \curvearrowright_{\beta} \Delta$ 
is an equivariantly approximable action by automorphisms, 
then the associated semidirect product $\Delta \rtimes_{\beta} \Gamma$ 
is a LEF group. Our next Proposition is a generalisation of their result. 

\begin{Thm}  \label{VerGorEAProp}
Let $\mathcal{C}$ be the class of LEF groups; let $\Gamma$ be a LEF group and suppose 
$\Gamma \curvearrowright_{\beta} \Delta$ is an $\mathcal{C}$-LEF action by automorphisms. 
Then $\Delta \rtimes_{\beta} \Gamma$ is a LEF group. 
\end{Thm}

\begin{proof}
Let $A \subseteq \Delta \rtimes_{\beta} \Gamma$ be a finite subset. 
Apply Lemma \ref{GenHayesSaleLem} to $A$ to obtain 
finite subsets $F_1 \subseteq \Delta$ and $F_2 \subseteq \Gamma$. 
Enlarging $F_1$ and $F_2$ (as we may, by Remark \ref{GenHayesSaleBiggerRmrk}), 
we may assume $A \subseteq F_1 \cdot F_2$ and $e \in F_2$. 
Let $Q$ be a finite group; $K$ be a LEF group; 
$Q \curvearrowright_{\gamma} K$ be an action by automorphisms, 
and $\pi: F_1 \hookrightarrow K$, $\rho: F_2 \rightarrow Q$ 
by partial homomorphisms satisfying the conditions of Definition \ref{LEFActGrpDef}. 
Let $d$ be the discrete metric on $G = K \rtimes_{\gamma} Q$, 
given by $d(g,h) = \delta_{g,h}$. 
Define $\sigma : \Delta \rtimes_{\beta} \Gamma \rightarrow G$ 
by $\sigma (h,g) = (\pi(h),\rho(g))$ for $(h,g) \in F_1 \cdot F_2$, 
and extended arbitrarily to $\Delta \rtimes_{\beta} \Gamma$. 
By construction, $\sigma$ satisfies the hypotheses of Lemma \ref{GenHayesSaleLem} 
with $\epsilon = 1/2$, so $\sigma |_A$ is a partial homomorphism. 
Further, the restriction of $\sigma$ to $F_1$ is injective, since $\pi$ is. 

Next, let $P$ be a finite group and let $\phi : F_2 \rightarrow P$ 
be an injective partial homomorphism. 
Then there is an injective partial homomorphism $\Phi : A \rightarrow G \times P$ 
given by $\Phi(h,g) = (\sigma(h,g),\phi(g)$. 
Note that $K$ is a finite-index subgroup of $G \times P$, so that by 
Proposition \ref{VLEFisLEFProp}, $G \times P$ is a LEF group. 
Since this holds for every such $A$, 
$\Delta \rtimes_{\beta} \Gamma$ is locally embeddable into $\mathcal{C}$. 
The conclusion now follows from Remark \ref{LELEFRem}. 
\end{proof}

Our main source of LEF actions on groups shall be the next Proposition. 

\begin{Prop} \label{LEFActEAAutoProp}
Let $\mathcal{L}$ be a set-theoretic
(respectively graph-theoretic) halo of groups and let $\mathcal{C}$ be a class of groups. 
Suppose that for every finite set (graph) $Y$, $L(Y)$ lies in $\mathcal{C}$. 
Let $X$ be a set (graph) and suppose $\Gamma \curvearrowright_{\alpha} X$ is a LEF action. 
Then the induced action $\hat{\alpha} : \Gamma \rightarrow \Aut(L(X))$
is a $\mathcal{C}$-LEF action. 
\end{Prop}

\begin{proof}
Let $F \subseteq \Gamma$ and $E \subseteq L(X)$ be finite sets.
Let $Z \subseteq X$ be a finite subset (induced subgraph)
such that $E \subseteq L(Z) \leq L(X)$.
Let $\rho : \Gamma \rightarrow Q$; 
$\pi : Z \hookrightarrow Y$ and
$Q \curvearrowright_{\beta} Y$ satisfy the conditions of
Definition \ref{LEFActionDef}. 
Then for all $g \in F$, the restrictions of $\pi \circ \alpha(g)$
and $\beta(g) \circ \pi$ to $Z \cap \alpha(g)^{-1} (Z)$ are equal.
Now for $h \in E$, if $\hat{\alpha}(g) (h) \in E$,
then since $E \subseteq L(Z)$ we have
$h \in L(Z) \cap \hat{\alpha}(g)^{-1} (L(Z)) = L (Z \cap \alpha(g)^{-1}(Z))$.
By functoriality, $\hat{\pi} \circ \hat{\alpha} (g)$
and $\hat{\beta} (g) \circ \hat{\pi}$ agree on
$L (Z \cap \alpha(g)^{-1}(Z))$,
so in particular
$$\hat{\pi} (\hat{\alpha} (g)(h))=\hat{\beta}(g)(\hat{\pi}(h)),$$
where $Q \curvearrowright_{\hat{\beta}} L(Y)$ 
is the induced action by automorphisms. 

Thus $\rho$; $\hat{\pi} : L(Z) \hookrightarrow L(Y)$
and $Q \curvearrowright_{\hat{\beta}} L(Y)$
satisfy Definition \ref{LEFActGrpDef}. 
\end{proof}

As an immediate consequence,
we may conclude that many halo products are LEF.

\begin{Cor} \label{LEFHaloProdCor}
Let $\mathcal{L}$ be a set-theoretic
(respectively graph-theoretic) halo of groups.
Let $X$ be a set (graph) and suppose that $L(X)$ is a LEF group.
Then for any LEF action $\Gamma \curvearrowright_{\alpha} X$,
the halo product $L(X) \rtimes_{\hat{\alpha}} \Gamma$
is a LEF group.
\end{Cor}

\begin{proof}
This follows from Proposition \ref{LEFActEAAutoProp} 
(applied with $\mathcal{C}$ the class of LEF groups) 
and Theorem \ref{VerGorEAProp}. 
\end{proof}

\begin{Cor}
Let $\mathcal{L}$ be a set-theoretic halo of groups.
Let $X$ be a set and suppose that $L(X)$ is a LEF group. 
Then for any action of a LERF group $\Gamma \curvearrowright_{\alpha} X$, 
the halo product $L(X) \rtimes_{\hat{\alpha}} \Gamma$ is a LEF group.
\end{Cor}

\begin{proof}
    This is immediate from Corollary \ref{LEFHaloProdCor} and 
    Theorem \ref{LERFLEFActThm}. 
\end{proof}

\subsection{Sofic actions on sets and graphs}

The concepts studied in this Subsection were introduced in 
\cite{GaoElayPatchSet} and \cite{GaoElayPatchGraph}. 

\begin{Def} \label{SoficActionDefn}
Let $\Gamma \curvearrowright_{\alpha} X$ be an action of a group on a set.
Let $A$ be a finite set and let $\varphi : \Gamma \rightarrow \Sym(A)$
be a function.
For finite subsets $F \subseteq \Gamma$, $E \subseteq X$ and $\epsilon >0$,
we say that $\varphi$ is an \emph{$(F,E,\epsilon)$-orbit approximation
of $\alpha$} if there exists a finite set $B$,
a subset $S \subseteq A$ with
$\lvert S \rvert > (1-\epsilon) \lvert A \rvert$,
and injective functions $(\pi_s : E \hookrightarrow B)_{s \in S}$ such that,
for all $g \in F$, $x \in E$ and $s \in S$,
if $\varphi(g) \cdot s \in S$ and $\alpha(g)^{-1} \cdot x \in E$ then:
$$\pi_{\varphi(g) \cdot s} (x) = \pi_s (\alpha(g)^{-1} \cdot x). $$
We describe $\alpha$ as a \emph{sofic action} if,
for all finite subsets $F \subseteq \Gamma$, $E \subseteq X$ and $\epsilon >0$,
there exists a finite set $A$ and a function
$\varphi : \Gamma \rightarrow \Sym(A)$ which is a unital
$(F,\epsilon)$-multiplicative $(F,E,\epsilon)$-orbit approximation
of $\alpha$.

Soficity of an action of a group on a graph is defined 
completely analogously, by requiring in the definition of 
an $(F,E,\epsilon)$-orbit approximation that $E$ is an induced subgraph of $X$; 
$B$ is a finite graph, and the $\pi_s$ are inclusions of induced subgraphs
\end{Def}

\begin{Rem} \label{SoficActBasicPropsRem}
    Sofic actions on sets (graphs) satisfy a suite of basic closure properties 
    closely analogous to those described for LEF actions in Remark \ref{LEFActionBasicsRem}. 
    They behave well with respect to: 
    pullbacks under homomorphisms; restriction to invariant subsets (induced subgraphs); 
    ascending chains of subgroups; actions on disjoint unions of sets; 
    graph-complements; passing from a graph to its vertex-set and passing from 
    a set to the complete or empty graph on that set 
    (see \cite{GaoElayPatchSet} Propositions 2.14 and 2.16, 
    and \cite{GaoElayPatchGraph} Proposition 1.2). 
\end{Rem}

The next result is Theorem 1.7 of \cite{GaoElayPatchGraph}. 

\begin{Thm} \label{SoficGrpsSoficActAmenStabThm}
    If $X = (V,E)$ is a graph; $\Gamma$ is a sofic group, 
    and $\Gamma \curvearrowright_{\alpha} X$ is a vertex-transitive action by 
    graph automorphisms such that the stabiliser of a vertex is amenable, 
    then $\alpha$ is a sofic action. 
    In particular (by Remark \ref{SoficActBasicPropsRem}), 
    any action of a sofic group on a set, with all point-stabilisers 
    being amenable subgroups, is a sofic action. 
\end{Thm}

\begin{Lemma} \label{RefinementLem}
Let $\varphi : \Gamma \rightarrow \Sym(A)$ be
a function, let $F \subseteq \Gamma$ be finite
and let $S \subseteq A$ with:
$$\lvert S \rvert > (1-\frac{\epsilon}{\lvert F \rvert+1}) \lvert A \rvert.$$
Then there exists $S_0 \subseteq S$ with $\lvert S_0 \rvert > (1-\epsilon) \lvert A \rvert$, such that for all $s \in S_0$ and $g \in F$,
we have $\varphi(g) \cdot s \in S$.
\end{Lemma}

\begin{proof}
We may take:
$$S_0 = S \cap \bigcap_{g \in F} \varphi(g)^{-1} (S)$$
so that $\varphi(g) \cdot s \in S$ for all $s \in S_0$ and $g \in F$,
and:
\begin{align*}
\big\lvert A \setminus \big(S \cap \bigcap_{g \in F} \varphi(g)^{-1} (S)\big) \big\rvert & = \big\lvert (A \setminus S) \cup \bigcup_{g\in F} (A \setminus \varphi(g)^{-1} (S) \big\rvert \\
& \leq (\lvert F \rvert + 1) \lvert A \setminus S \rvert \\
& < \epsilon \lvert A \rvert
\end{align*}
as required.
\end{proof}

\begin{Prop} \label{LEFActSoficProp}
Suppose $\Gamma \curvearrowright_{\alpha} X$ is a LEF action.
Then for any finite subsets $F \subseteq \Gamma$
and $Z \subseteq X$,
there exist finite sets $A$ and $B$; injections $(\pi_a : Z \hookrightarrow B)_{a \in A}$
and a unital, $(F,\epsilon)$-multiplicative function
$\varphi : \Gamma \rightarrow \Sym(A)$ such that
satisfying the conditions of Definition \ref{SoficActionDefn}
for all $\epsilon > 0$, with $S=A$.
In particular, $\Gamma \curvearrowright_{\alpha} X$ is a sofic action.
\end{Prop}

\begin{proof}
Suppose, as we may, that $e \in F$. 
Let $Q$ be a finite group; $Y$ be a finite set; $Q \curvearrowright_{\beta} Y$ 
and $\rho: \Gamma \rightarrow Q$, $\pi : Z \hookrightarrow Y$ 
be as in Definition \ref{LEFActionDef}, 
so that there is an induced homomorphism $\hat{\beta} : Q \rightarrow \Sym(Y)$. 
Set $A=\Sym(Y)$ and let $\lambda : A \rightarrow \Sym(A)$ be the left-regular permutation 
representation. 
Let $\varphi = \lambda \circ \hat{\beta} \circ \rho : \Gamma \rightarrow \Sym(A)$. 
Then $\varphi$ is unital and 
$(F,\epsilon)$-multiplicative, for all $\epsilon > 0$. 
Set $B=Y$ and for $a \in A$ let $\pi_a : Z \hookrightarrow B$ be given by 
$\pi_a (x) = a^{-1} (\pi(s))$. 
Then for $x \in Z$; $g \in F$ and $a \in A$, if $\alpha (g)^{-1} (x) \in Z$, we have: 
$$\pi_a (\alpha(g)^{-1}(x)) = (a^{-1} \hat{\beta}(\rho(g))^{-1})(\pi(x)) = (\hat{\beta}(\rho(g)))^{-1}(\pi(x)) = \pi_{\varphi (g) (a)} (x)$$
since $\varphi (g) (a) = \lambda(\hat{\beta}(\rho(g)))[a] = \hat{\beta}(\rho(g)) a$. 
\end{proof}

\begin{Rem}
    Combining Proposition \ref{LEFActSoficProp} with Theorem \ref{LERFLEFActThm}, 
    we recover the conclusion that every action of a LERF group on a set is sofic: 
    this was the main result of \cite{Gao}. 
\end{Rem}

\begin{Ex}
    Not every sofic action on a set is an LEF action. 
    Let $\Gamma$ be a finitely presented amenable group which is not residually finite. 
    Then by Proposition \ref{fpLEFRFProp}, $\Gamma$ is not an LEF group, 
    so by Proposition \ref{LEFGrpFaithLEFActProp}, $\Gamma$ admits no faithful LEF action. 
    On the other hand, the action of $\Gamma$ on itself by left-multiplication 
    is a sofic action, by Theorem \ref{SoficGrpsSoficActAmenStabThm}. 
    An example of such a group $\Gamma$ is described in Subsection 5.4 
    of \cite{CornGuyoPits}: the group $A_4 / Z$ constructed therein is 
    finitely presented; $3$-step soluble (hence amenable) and non-Hopfian 
    (hence not residually finite, as it is finitely generated). 
\end{Ex}

\subsection{Sofic actions on groups}

We introduce the following notion of a sofic action of one group on another by automorphisms, 
which is closely inspired by Definition \ref{SoficActionDefn}, 
and bears much the same relation to it that Definition \ref{LEFActGrpDef} 
did to Definition \ref{LEFActionDef}. 

\begin{Def} \label{SoficAutoActionDefn}
Let $\Gamma , \Delta$ be groups, 
let $\alpha : \Gamma \rightarrow \Aut (\Delta)$ 
be a homomorphism and let $\varphi : \Gamma \rightarrow \Sym(A)$ 
be a function. 
For $F \subseteq \Gamma$ and $E \subseteq \Delta$ finite subsets; 
$\epsilon > 0$ and $\Lambda$ a group, 
we say that $\varphi$ is an $(F,E,\epsilon,\Lambda)$-automorphic 
approximation of $\alpha$ if there exists $S \subseteq A$ 
with $\lvert S \rvert > (1-\epsilon)\lvert A \rvert$ 
and injective partial homomorphisms 
$(\pi_s : E \hookrightarrow \Lambda)_{s \in S}$ 
such that for all $g \in F$, $h \in E$ and $s \in S$, 
if $\varphi (g) \cdot s \in S$ and $\alpha(g)^{-1} [h] \in E$, then: 
$$\pi_{\varphi (g) \cdot s} (h) = \pi_s (\alpha(g)^{-1} [h]).$$
For $\mathcal{C}$ a class of groups, 
we say that $\alpha$ is a sofic $\mathcal{C}$-action if, 
for all finite subsets $F \subseteq \Gamma$, $E \subseteq \Delta$ 
and $\epsilon > 0$, there exists a group $\Lambda \in \mathcal{C}$, 
a finite set $A$ and a unital $(F,\epsilon)$-multiplicative 
$(F,E,\epsilon,\Lambda)$-automorphic approximation of $\alpha$. 
\end{Def} 



\begin{Prop} \label{AmenSoficAutProp}
Let $\Gamma$ be an amenable group. 
Then for any group $\Delta$ and any class $\mathcal{C}$ of groups for which 
$\Delta$ is locally embeddable into $\mathcal{C}$, 
every action $\Gamma \curvearrowright_{\alpha} \Delta$ 
by automorphisms is a sofic $\mathcal{C}$-action. 
\end{Prop}

\begin{proof}
Let $F \subseteq \Gamma$ and $E \subseteq \Delta$ be nonempty 
finite subsets and let $\epsilon > 0$. 
Since $\Gamma$ is amenable, there exists a finite subset $A \subseteq \Gamma$ 
such that for all $g \in F$, 
$\lvert g A \setminus A \rvert < \epsilon \lvert A \rvert / 2 \lvert F \rvert$. 
For $g \in F$ and $a \in A \cap g^{-1} A$, set $\varphi(g)[a] = ga \in A$. 
Then $\varphi(g)$ is a partially-defined permutation of $A$, 
so we may extend it (arbitrarily) to an element of $\Sym(A)$, 
also denoted $\varphi(g)$. 
Extending $\varphi$ arbitrarily to $\Gamma$, we have that $\varphi : \Gamma \rightarrow \Sym(A)$ 
is $(F,\epsilon)$-multiplicative. 
Moreover, letting $S \subseteq A$ be given by 
$$S = A \cap \bigcap_{g \in F} g^{-1} A,$$
we have $\lvert S \rvert > (1-\epsilon) \lvert A \rvert$ and 
for all $s \in S$ and $g \in F$, $\varphi(g)[s] = gs$. 

For $s \in S$ define $\tilde{\pi}_s : E \hookrightarrow \Delta$ by 
$\tilde{\pi}_s [h] = \alpha(s)^{-1}(h)$. 
Then $\tilde{\pi}_s$ is an injective partial homomorphism 
and for all $g \in F$, $h \in E$ and $s \in S$ for which $\varphi(g)[s] \in S$ and 
$\alpha(g)^{-1}[h] \in E$, we have: 
$$\tilde{\pi}_{\varphi(g)[s]} (h) = \tilde{\pi}_{gs} (h) = \alpha(gs)^{-1}[h] = \alpha(s)^{-1}[\alpha(g)^{-1}[h]] = \tilde{\pi}_s (\alpha(g)^{-1}[h]).$$
Finally, let $B  = \bigcup_{s \in S} \alpha(s)^{-1}(E)$, a finite subset of 
$\Delta$, and let $\psi : B \hookrightarrow \Lambda$ be an injective 
partial homomorphism for some group $\Lambda \in \mathcal{C}$. 
Then $(\pi_s = \psi \circ \tilde{\pi}_s : E \hookrightarrow \Lambda)_{s \in S}$ 
are partial homomorphisms and witness that $\varphi$ 
is an $(F,E,\epsilon,\Lambda)$-automorphic approximation of $\alpha$. 
\end{proof}

\begin{Prop} \label{EAActionsSoficProp}
Let $\mathcal{C}$ be a class of groups and let $\Gamma, \Delta$ be groups. 
If $\Gamma \curvearrowright_{\alpha} \Delta$ is a $\mathcal{C}$-LEF action, 
then $\alpha$ is a sofic $\mathcal{C}$-action. 
\end{Prop}

\begin{proof}
Let $F \subseteq \Gamma$ and $E \subseteq \Delta$ be finite subsets. 
Enlarging $F$, as we may, we assume $e \in F = F^{-1}$. 
Let $Q$ and $K$ be finite groups; let $\gamma : Q \rightarrow \Aut(K)$ 
be an action of $Q$ on $K$ by group automorphisms, 
and $\rho : F \rightarrow Q$ and $\pi : E \hookrightarrow K$ 
be as in Definition \ref{LEFActGrpDef}. 
Set $A = \Aut(K)$, and let $\lambda : A \rightarrow \Sym(A)$ 
be the left-regular permutation representation 
and take $\varphi : \Gamma \rightarrow \Sym(A)$ 
to be an arbitrary extension of $\lambda \circ \gamma \circ \rho$ to $\Gamma$. 
Set $B=K$ and for $a\in A$, let $\pi_a : E \hookrightarrow K$ 
be given by $\pi_a (h) = a^{-1}(\pi(h))$. Then 
for $g \in F$, $h \in E$ and $a \in A$, provided $\alpha(g)^{-1}[h] \in E$ we have: 
\begin{align*}
\pi_a (\alpha(g)^{-1}[h])=a^{-1} (\pi(\alpha(g^{-1})[h]))&=a^{-1} \gamma (\rho(g))^{-1}[\pi(h)] \\
&= (\gamma (\rho(g)) a)^{-1} [\pi(h)]\\ 
&= (\varphi(g)[a])^{-1} [\pi(h)] \\
&= \pi_{\varphi(g)[a]} (h).
\end{align*}
Thus $\varphi$ is an $(F,E,\epsilon,K)$-automorphic approximation of $\alpha$ 
for all $\epsilon > 0$. Since $F$ and $E$ were arbitrary, we have the desired conclusion. 
\end{proof}

The key source of sofic actions of one group on another 
for the next Section shall be the following Proposition. 

\begin{Prop} \label{SoficActSoficAutoProp}
Let $\mathcal{C}$ be a class of groups 
and let $\mathcal{L}$ be a set-theoretic 
(respectively graph-theoretic) halo of groups. 
Suppose that for every finite set 
(graph) $X$ the group $L(X)$ lies in $\mathcal{C}$. 
Then for any set (graph) $X$ and and sofic group action 
$\Gamma \curvearrowright_{\alpha} X$, 
the induced action $\hat{\alpha} : \Gamma \rightarrow \Aut(L(X))$ 
is a sofic $\mathcal{C}$-action. 
\end{Prop}

\begin{proof}
Let $F \subseteq \Gamma$ and $E \subseteq L(X)$ be finite sets 
and let $\epsilon > 0$. 
Let $Y \subseteq X$ be a finite subset (induced subgraph) 
such that $E \subseteq L(Y) \leq L(X)$; 
let $A$ be a finite set and let $\varphi:\Gamma\rightarrow\Sym(A)$ 
be a unital $(F,\epsilon)$-multiplicative 
$(F,E,\epsilon)$-orbit approximation, with associated data 
$S \subseteq A$ and $(\pi_s : Y \hookrightarrow B)_{s \in S}$ 
for some finite set (graph) $B$, satisfying the conditions of 
Definition \ref{SoficActionDefn}. 

We have an induced action $\Gamma \curvearrowright_{\hat{\alpha}}$ 
by automorphisms and injective homomorphisms 
$\hat{\pi}_s : L(Y) \hookrightarrow L(B)$. 
For $g \in F$, $h \in E$ and $s \in S$, 
if $\varphi(g) \cdot s \in S$ and $\hat{\alpha} (g)^{-1} (h) \in E$ 
we have 
$h \in L(Y) \cap \hat{\alpha}(g)(L(Y)) = L(Y \cap \alpha(g)(Y))$. 
The restriction of $\pi_{\varphi(g)\cdot s}$ 
to $Y \cap \alpha(g)(Y)$ agrees with $\pi_s \circ \alpha(g)^{-1}$ 
so by functoriality, 
$\hat{\pi}_{\varphi(g)\cdot s}$ agrees with 
$\hat{\pi}_s \circ \hat{\alpha}(g)^{-1}$ 
on $L(Y \cap \alpha(g)(Y)) \ni h$, 
that is $\hat{\pi}_{\varphi(g)\cdot s} (h) 
= \hat{\pi}_s ( \hat{\alpha}(g)^{-1}(h))$. 
Since $L(B)$ lies in $\mathcal{C}$, 
we have that $\hat{\alpha}$ satisfies the conditions of 
Definition \ref{SoficAutoActionDefn}. 
\end{proof}

\section{Constructing metric approximations of semidirect products}

\begin{Thm} \label{BigTechThm}
Let $\mathcal{G} = (G_i,d_i)_{i \in I}$ be a metric family
of groups.
Suppose $\mathcal{G}$ is product-compatible and wreath-compatible. 
Let $\mathcal{C}$ be the class of $\mathcal{G}$-approximable groups. 
For $c > 0$ let $\mathcal{C}_c$ be the class of $(\mathcal{G},c)$-approximable groups. 
Let $\Gamma , \Delta \in \mathcal{C}$,
and let $ \Gamma \curvearrowright_{\beta} \Delta$
be a sofic $\mathcal{C}_{c}$-action by automorphisms for some $c>0$.
Then $\Delta \rtimes_{\beta} \Gamma \in \mathcal{C}$.
\end{Thm}

Before proceeding with the proof of Theorem \ref{BigTechThm}
we record some consequences.

\begin{Cor} \label{BigSemiCoroll}
Let $\mathcal{C}$ be one of the following:
\begin{itemize}
\item[(a)] The class of sofic groups;

\item[(b)] The class of $\mathbb{K}$-linear sofic groups, for $\mathbb{K}$ a field; 

\item[(c)] The class of hyperlinear groups;

\item[(d)] The class of weakly sofic groups. 

\end{itemize}
Let $\Gamma , \Delta \in \mathcal{C}$,
and let $\beta : \Gamma \curvearrowright \Delta$
be a sofic $\mathcal{C}$-action by automorphisms.
Then $\Delta \rtimes_{\beta} \Gamma \in \mathcal{C}$.
\end{Cor}

\begin{proof}
Let $\mathcal{G}$ be, respectively, 
$\mathcal{G}_{\sof}$; $\mathcal{G}_{\mathbb{K}}$; $\mathcal{G}_{\hyp}$ 
or $\mathcal{G}_{\weak}$. 
Then by Lemma \ref{SoficUnifSepLemma}; Lemma \ref{LinSoficUnifSepLemma}; 
Lemma \ref{HypUnifSepLemma} or Lemma \ref{WSoficUnifSepLemma}, 
there exists $c>0$ such that every $\mathcal{G}$-approximable group 
is $(\mathcal{G},c)$-approximable. 
Thus by Theorem \ref{BigTechThm}, it is sufficient to check that
the relevant matric families are product- and wreath-compatible.
This follows from Lemma \ref{SoficCompatLem} in case (a);
Lemma \ref{LinSoficCompatLemma} in case (b); Lemma \ref{HypCompatLem} in case (c) 
and Lemma \ref{WSoficCompatLem} in case (d).
\end{proof}

By Proposition \ref{AmenSoficAutProp}, 
we immediately deduce the following from Corollary \ref{BigSemiCoroll}, 
which recovers special cases of results from 
\cite{ElekSzaboSofic} (for sofic groups); \cite{ArzhPaunLinear,Stolz} (linear sofic groups); 
\cite{ArzhBerlFiSeGlen,HoltRees} (hyperlinear groups), and \cite{HoltRees} (weakly sofic groups): 
the cited works show that any $\mathcal{C}$-by-amenable group lies in $\mathcal{C}$ 
for $\mathcal{C}$ one of the classes (a)-(d) from Corollary \ref{BigSemiCoroll}. 

\begin{Cor} \label{AmenQuotProp}
    Let $\mathcal{C}$ be 
    Let $\Gamma$ be an amenable group and let $\Delta \in \mathcal{C}$. 
    Then for any homomorphism $\beta : \Gamma \rightarrow \Delta$, 
    we have $\Delta \rtimes_{\beta} \Gamma \in \mathcal{C}$.
\end{Cor}

\begin{Cor} \label{BigHaloCorollGen}
Let $\mathcal{G}$; $\mathcal{C}$ and $\mathcal{C}_c$ be as in Theorem \ref{BigTechThm};
let $\mathcal{L}$ be a set-theoretic
(respectively graph-theoretic) halo of groups,
and suppose there exists $c>0$ such that for any finite set (graph) $B$,
$L(B) \in \mathcal{C}_c$.
Then for any set (graph) $X$, any $\Gamma \in \mathcal{C}$,
and any sofic action $\Gamma \curvearrowright_{\alpha} X$,
the associated halo product $L(X) \rtimes_{\hat{\alpha}} \Gamma$
lies in $\mathcal{C}$.
\end{Cor}

\begin{proof}
This follows from Theorem \ref{BigTechThm} and
Proposition \ref{SoficActSoficAutoProp} (the latter applied with $\mathcal{C}_c$ 
in place of $\mathcal{C}$).
\end{proof}

\begin{Cor} \label{BigHaloCorollSpec}
Let $\mathcal{C}$ be one of the classes (a)-(d) in Corollary 
\ref{BigSemiCoroll}, and let $\mathcal{L}$ be as in Corollary
\ref{BigHaloCorollGen}.
Then for any set (graph) $X$, any $\Gamma \in \mathcal{C}$,
and any sofic action $\Gamma \curvearrowright_{\alpha} X$,
the associated halo product $L(X) \rtimes_{\hat{\alpha}} \Gamma$
lies in $\mathcal{C}$.
\end{Cor}

\begin{proof}
We deduce this from Corollary \ref{BigHaloCorollGen}
just as Corollary \ref{BigSemiCoroll} was deduced from
Theorem \ref{BigTechThm}.
\end{proof}

\begin{proof}[Proof of Theorem \ref{BigTechThm}]
Let $F \subseteq \Delta \rtimes_{\beta} \Gamma$ be a finite subset
and let $\epsilon > 0$.
Let $F_1 \subseteq \Delta$, $F_2 \subseteq \Gamma$ be finite
subsets as in Lemma \ref{GenHayesSaleLem} (Hayes and Sale).
As per Remark \ref{GenHayesSaleBiggerRmrk}, we may enlarge $F_1$ and $F_2$ as required.
In particular, we may assume that $e \in F_1$; 
$e \in F_2 = F_2 ^{-1}$
and $F \subseteq F_1 \cdot F_2$.
Let $E = \lbrace \beta(g)[h] : g \in F_2 , h \in F_1 \rbrace$. 
Since $\Gamma$ is $\mathcal{G}$-approximable, 
there exists $c' > 0$ such that for every $F' \subseteq \Gamma$ finite 
and every $\epsilon' > 0$, there exists an $(F',\epsilon',c')$-$\mathcal{G}$-representation 
of $\Gamma$. 
Let $\mu > 0$ be such that, for every $i,j \in I$ 
there exists $k \in I$ and a $\mu$-product map $G_i \times G_j \rightarrow G_k$. 

We shall in the course of the proof introduce, in order, quantities 
$\epsilon_1 , \epsilon_2 , \epsilon_3 , \epsilon_4 > 0$, 
where $\epsilon_1$ is arbitrary and 
each $\epsilon_{i+1}$ is a function of $\epsilon_i$ 
and the quantities fixed in the previous paragraph 
(chosen in the course of the proof), 
such that $\epsilon_{i+1}$ can be made arbitrarily small 
by choosing $\epsilon_i$ sufficiently small. 
We show that there exists a constant $c'' > 0$ 
such that there exists an $(F,\epsilon_4,c'')$-$\mathcal{G}$-representation 
of $\Delta \rtimes_{\beta} \Gamma$. 
Thus, choosing $\epsilon_1$ sufficiently small, such that $\epsilon_4 < \epsilon$, 
we have an $(F,\epsilon,c'')$-$\mathcal{G}$-representation. 
Since $F$ and $\epsilon$ were arbitrary, 
we conclude that $\Delta \rtimes_{\beta} \Gamma$ is $(\mathcal{G},c'')$-approximable. 
In fact, in our proof it shall be possible to take $c'' = \mu \min (c/2,c')$. 

For some $\Lambda \in \mathcal{C}_c$
we have $\varphi : \Gamma \rightarrow \Sym(A)$, 
a unital $(F_2,\epsilon_1)$-multiplicative
$(F_2,E,\epsilon_1,\Lambda)$-automorphic approximation
to $\beta$, with associated data $S \subseteq A$ and
injective partial homomorphisms
$(\pi_s : E \hookrightarrow \Lambda)_{s\in S}$.
By Lemma \ref{RefinementLem}, 
we have a subset $S_0 \subseteq S$ such that
$\varphi(g)\cdot s$ for all $g \in F_2$ and $s \in S_0$,
and $\lvert S_0 \rvert > (1-\epsilon_2) \lvert A \rvert$.
In particular, 
\begin{equation}
\pi_{\varphi(g)\cdot s} (h) = \pi_s (\beta(g)^{-1}[h])
\end{equation}
for all $s \in S_0$; and $g \in F_2$ and $h \in E$ 
for which $\beta(g)^{-1}[h] \in E$. 

Set: $$F_3 = \bigcup_{s \in S_0} \pi_s (E),$$ 
a finite subset of $\Lambda$ 
and let $\sigma_E : \Lambda \rightarrow G_{i_1}$ 
be a $(F_3,\epsilon_1,c)$-$\mathcal{G}$-representation. 
Define $\hat{\sigma}_E : E \rightarrow G_{i_1} ^A$ by: 
$$\hat{\sigma}_E (h)_a = \left\{ \begin{array}{ll}
\sigma_E (\pi_a (h)) & a \in S_0 \\
1_{G_{i_1}} &  a \in A \setminus S_0 
\end{array} \right.$$
Let $\tau : G_{i_1} ^A \rightarrow G_{i_2}$ 
and $\psi : \Sym(A) \rightarrow G_{i_2}$ be, respectively, 
a base map and an acting map 
as in the definition of wreath compatibility, 
and define $\Phi : \Delta \rtimes_{\beta} \Gamma \rightarrow G_{i_2}$ by: 
$$\Phi (h,g) = \tau (\hat{\sigma}_E (h)) \psi (\varphi(g)).$$
First, we note that $\Phi$ is unital, since 
$\tau$, $\hat{\sigma}_E$, $\psi$ and $\varphi$ are. 

Second, we verify that $\Phi$ is 
$(F,\epsilon_3)$-multiplicative, 
by verifying that $F_1$ and $F_2$ satisfy the conditions of Lemma 
\ref{GenHayesSaleLem}. 

For condition (i), if $h_1 , h_2 \in F_1$ are such that 
$h_1 h_2 \in F_1$, then for $a \in S_0$, 
we have $\pi_a (h_1) , \pi_a (h_2) , \pi_a (h_1 h_2) \in F_3$ 
and $\pi_a (h_1 h_2) = \pi_a (h_1) \pi_a (h_2)$, 
so: 
$$d_{i_1} (\sigma_E (\pi_a (h_1 h_2)),\sigma_E (\pi_a (h_1))\sigma_E (\pi_a (h_2))) < \epsilon_1.$$
Set $X = S_0$ in part (a) of the definition of wreath-compatibility, so that: 
$$d_{i_2} (\Phi(h_1 h_2 , e) , \Phi(h_1, e)\Phi(h_2 , e)) 
= d_{i_2} (\tau(\hat{\sigma}_E (h_1 h_2)),\tau(\hat{\sigma}_E (h_1))\tau(\hat{\sigma}_E (h_2)))$$
For condition (ii), for $g_1 , g_2\in F_2$ with $g_1 g_2\in F_2$, 
$$d_A (\varphi(g_1 g_2) , \varphi(g_1)\varphi(g_2) ) < \epsilon_1$$
Then by condition (b) of wreath-compatibility, 
$$d_{i_2} (\Phi(e,g_1g_2),\Phi(e,g_1)\Phi(e,g_2)) 
= d_{i_2} (\psi(\varphi(g_1 g_2)),\psi(\varphi(g_1 ))\psi(\varphi(g_2))) < \epsilon_2.$$
For condition (iii), 
$\Phi(h,e) = \tau(\hat{\sigma}_E(h))$ 
and $\Phi(e,g) = \psi(\varphi(g))$ 
so $\Phi(h,e)\Phi(e,g) = \Phi(h,g)$. 

For condition (iv), for $h \in F_1$ and $g \in F_2$ 
we are comparing: 
$$\Phi(e,g)\Phi(h,e) = \psi(\varphi(g)) \tau (\hat{\sigma}_E(h))$$
with: 
$$\Phi(\beta(g)[h],e) \Phi(e,g) 
= \tau(\hat{\sigma}_E (\beta(g)[h])) \psi(\varphi(g)),$$
that is, writing $\hat{\sigma}_E (\beta(g)[h])_a =k_a\in G_{i_1}$, 
we are comparing $\tau(\hat{\sigma}_E (h))$ with: 
$$\psi(\varphi(g))^{-1} \tau(\hat{\sigma}_E (\beta(g)[h])) \psi(\varphi(g)) = \tau((k_{\varphi(g)[a]})_a)$$
(the latter equality holding 
by condition (c) of wreath-compatibility); 
we seek to prove they are $\epsilon_3$-close 
in the metric $d_{i_2}$. 
Now $\lvert S_0 \rvert > (1-\epsilon_2) \lvert A \rvert$ 
so: $$\lvert S_0 \cap \varphi(g)^{-1}(S_0) \rvert > (1-2\epsilon_2) \lvert A \rvert$$
hence by condition (a) of wreath-compatibility, 
it suffices to check that for all
$a \in S_0 \cap \varphi(g)^{-1}(S_0)$, 
we have $\hat{\sigma}_E (h)_a = k_{\varphi(g)[a]}$, 
that is, that we have: 
$$\sigma_E (\pi_a(h)) 
= \hat{\sigma}_E (\beta(g)[h])_{\varphi(g)[a]}
= \sigma_E (\pi_{\varphi(g)[a]} (\beta(g)[h])).$$
But since $h \in F_1$, $\beta(g)[h] \in E$, 
by Definition \ref{SoficAutoActionDefn} we have: 
$$\pi_{\varphi(g)[a]} (\beta(g)[h])
= \pi_a(\beta(g)^{-1}[\beta(g)[h]])
= \pi_a (h)$$
as desired. 

Third, we claim that $\Phi$ is $(F_1,(c/2))$-separating 
(recalling that $F_1 \subseteq F$). 
Let $e \neq h \in F_1$. For $a \in S_0$ we have: 
$$d_{i_1} (\hat{\sigma}_E(h)_a,1_{G_{i_1}}) = d_i (\sigma_E (\pi_A(h)),1_{G_{i_1}}) > c$$
since $e \neq \pi_a (h) \in F_3$ and $\sigma_E$ 
is $(F_3,c)$-separating. 
By condition (a) of wreath-compatibility, 
applied with $X=S_0$, since 
$\lvert X \rvert > (1-\epsilon_2)\lvert A \rvert$ we have: 
$$d_{i_2} (\tau(\hat{\sigma}_E (h)),e_{g_{i_2}}) > (1-\epsilon_3) c > c/2$$
provided $\epsilon_3$ is sufficiently small. 

Finally, let $\theta : \Gamma \rightarrow G_{i_3}$ be 
an $(F_2,\epsilon_3,c')$-representation, for some $i_3 \in I$. 
let $\tilde{\theta} : \Delta \rtimes_{\beta} \Gamma \rightarrow G_{i_3}$ 
be given by $\tilde{\theta} (h,g)=\theta(g)$. 
Then, since $F_2 \subseteq F \subseteq F_1 F_2$, 
$\tilde{\theta}$ is unital and $(F,\epsilon_3)$-multiplicative. 
Let $\Delta_{i_2 , i_3} : G_{i_2} \times G_{i_3} \rightarrow G_{i_4}$ 
be a $\mu$-product map, and let 
$\Psi : \Delta \rtimes_{\beta} \Gamma \rightarrow G_{i_4}$ 
be given by $\Psi (h,g) = \Delta_{i_2 , i_3}  (\Phi(h,g),\tilde{\theta}(h,g))$. 
We claim that $\Psi$ is an $(F,\epsilon_{p-3},c'')$-representation, for $c'' = \mu \min (c/2,c')$, 
which shall complete the proof. 
Certainly $\Psi$ is unital (since $ \Delta_{i_2 , i_3}$, $\Phi$ and $\theta$ are), 
and $(F,\epsilon_4)$-multiplicative by Lemma \ref{ProdMapAlmostHomLem}. 
For $e \neq (h,g) \in F$, we have $h \in F_1$ and $g \in F_2$. 
If $g \neq e$, then $d_{i_3} (\tilde{\theta}(h,g),e)=d_{i_3} (\theta(g),e) \geq c'$, 
so by item (b) of Definition \ref{ProdCompatDef}, 
$d_{i_4} (\Psi (h,g),e) \geq \mu c'$. 
On the other hand, if $g = e$, then $e \neq h \in F_1$, 
so $d_{i_2} (\Phi(h,g),e) = d_{i_2} (\Phi(h,e),e) \geq c/2$, 
so by Definition \ref{ProdCompatDef} (b) once again, 
$d_{i_4} (\Psi (h,g),e) \geq \mu c/2$, 
so in any case $(h,g)$ is $c''$-separated by $\Psi$, as desired. 
\end{proof}

\section{Examples} \label{ExamplesSect}

\subsection{Symmetric and alternating enrichments}

For $X$ a set, let $\Sym_f(X)$ be the set of finitely supported 
permutations of $X$, and let $\Alt_f(X) \leq \Sym_f(X)$ be 
the subgroup of even permutations in $\Sym_f(X)$. 
It is easy to see that the assignment $X \mapsto \Sym_f(X)$ 
yields a set-theoretic halo of groups $\mathcal{S}$, 
with the inclusion of sets $Y \hookrightarrow X$ 
inducing a monomorphism of groups $\Sym_f(Y) \hookrightarrow \Sym_f(X)$ 
via the trivial extension of a permutation from $Y$ to $X$. 
Likewise the assignment $X \mapsto \Alt_f(X)$ 
gives rise to a set-theoretic halo of groups $\mathcal{A}$. 
A group action $\Gamma \curvearrowright_{\alpha} X$ 
yields an action $\hat{\alpha}$ of $\Gamma$ on $\Sym_f(X)$ by automorphisms, 
in accordance with Definition \ref{FunctHaloProdDefn}, via: 
$$\hat{\alpha}(g) : \sigma \mapsto \alpha(g) \sigma \alpha(g) ^{-1}$$ 
(for $g \in \Gamma$ and $\sigma \in \Sym_f(X)$), 
under which action $\Alt_f(X)$ is preserved by $\Gamma$. 
We thus have associated semidirect products: 
$\mathcal{S}(\Gamma,X) = \Sym_f(X) \rtimes_{\hat{\alpha}} \Gamma$ 
and $\mathcal{A}(\Gamma,X) = \Alt_f(X) \rtimes_{\hat{\alpha}} \Gamma$ 
called the \emph{symmetric enrichment} and \emph{alternating enrichment} 
of the action $\Gamma \curvearrowright_{\alpha} X$. 
Note that $\mathcal{A}(\Gamma,X)$ is an index-two 
subgroup of $\mathcal{S}(\Gamma,X)$. 
We follow the terminology of \cite{BradLEFExtn}, 
which was proposed by J.S. Wilson. 
In \cite{GeneTess}, $\mathcal{S}(\Gamma,X)$ 
was called the \emph{lampshuffler group} 
associated to the action $\Gamma \curvearrowright_{\alpha} X$. 
The action of $\Gamma$ on $X$ also gives rise, 
for every $n \in \mathbb{N}$, 
to an action on $X_n = X \times \lbrace 1 , \ldots , n \rbrace$ 
(acting on $X$ by $\alpha$, and trivially on the second factor) 
and the group $\mathcal{S}(\Gamma,X_n)$ was 
referred to in \cite{GeneTess} as the \emph{lampjuggler group}. 

Thanks to Corollaries \ref{BigHaloCorollSpec} and \ref{LEFHaloProdCor}, 
we have the following results, the latter of which 
was already proved as part of Theorem 3.3 in \cite{Brad}. 

\begin{Thm} \label{SoficSymEnrThm}
Suppose that $\Gamma$ is a sofic 
(respectively $\mathbb{K}$-linear sofic; hyperlinear; weakly sofic) 
group and that 
$\Gamma \curvearrowright_{\alpha} X$ is a sofic action on the set $X$. 
Then $\mathcal{S}(\Gamma,X)$ and $\mathcal{A}(\Gamma,X)$ 
are sofic ($\mathbb{K}$-linear sofic; hyperlinear; weakly sofic) groups. 
\end{Thm}

\begin{Thm} \label{SymEnrLEFThm}
Suppose that $\Gamma$ is an LEF group and that 
$\Gamma \curvearrowright_{\alpha} X$ is a LEF action on the set $X$. 
Then $\mathcal{S}(\Gamma,X)$ and $\mathcal{A}(\Gamma,X)$ 
are LEF groups. 
\end{Thm}

We also note that it is easy to determine when 
$\mathcal{S}(\Gamma,X)$ and $\mathcal{A}(\Gamma,X)$ 
are residually finite. 

\begin{Prop} \label{SymEnrRFProp}
The group $\mathcal{S}(\Gamma,X)$ (respectively $\mathcal{A}(\Gamma,X)$) 
is residually finite 
iff $\Gamma$ is a residually finite group and $X$ is a finite set. 
\end{Prop}

\begin{proof}
If $X$ is infinite, then $\mathcal{A}(\Gamma,X) \leq \mathcal{S}(\Gamma,X)$ 
contain the infinite simple subgroup $\Alt_f(X)$, 
hence cannot be residually finite. 
If $\Gamma$ is not residually finite, 
then neither is $\mathcal{A}(\Gamma,X)$, 
since $\Gamma$ embeds as a subgroup of $\mathcal{A}(\Gamma,X)$. 

Conversely suppose $X$ is finite (so that $\Sym_f(X) = \Sym(X)$) and $\Gamma$ 
is residually finite. Let $(\sigma , g) \in \mathcal{S}(\Gamma,X)$ 
be nontrivial. If $g \neq e$, then $g \notin N$ for some 
$N \vartriangleleft_f \Gamma$. 
Thus $\Sym(X) N \vartriangleleft_f  \mathcal{S}(\Gamma,X)$ 
with $(\sigma , g) \notin \Sym(X) N$. 
If on the other hand $g = e$, so that $\sigma \neq \id_X$, 
let $K$ be the kernel of the action of $\Gamma$ on $\Sym(X)$, 
so that $K \vartriangleleft_f \Gamma$. 
Then $\Sym(X)\rtimes (K\Gamma)$ is a finite quotient of $\mathcal{S}(\Gamma,X)$ 
in which $(\sigma , e)$ survives. 
\end{proof}

\subsection{Linear enrichments} \label{LinEnrSubsect}

Let $R$ be a commutative unital ring and let $R[X]$ be the free $R$-module 
with basis the set $X$. 
For a subset $Y \subseteq X$, $R[Y]$ is naturally a submodule of $R[X]$. 
Let $\GL_f (X,R)$ denote the subgroup of the $R$-module 
automorphism group $\Aut_R (R[X])$ of $R[X]$ 
consisting of those automorphisms $A \in \Aut_R (R[X])$ 
with the property that there exists some finite $Z_A \subseteq X$ 
for which $A$ preserves $R[Z_A] \leq R[X]$ and fixes all basis vectors $x \in X \setminus Z_A$. 
Then the inclusion of the subset $Y \subseteq X$ induces 
an inclusion of $\GL_f (Y,R)$ as a subgroup of $\GL_f (X,R)$ 
(extending $A \in \GL_f (Y,R)$ to an element of $\GL_f (X,R)$ by declaring $A$ 
to fix all basis vectors in $X \setminus Y$). 
The assignment $X \mapsto \GL_f (X,R)$ then clearly gives rise to a set-theoretic halo of groups. 

A group action $\Gamma \curvearrowright_{\alpha} X$ on the set $X$ 
induces an action of $\Gamma$ on $R[X]$ by permuting basis elements 
(also denoted by $\alpha$), 
and thus an action $\hat{\alpha}$ on $\Aut_R (R[X])$ given by: 
$$\hat{\alpha}(g) : A \mapsto \alpha (g) A \alpha (g) ^{-1}$$
(for $g \in \Gamma$ and $A \in \Aut_R (R[X])$). 
It is clear that $\GL_f (X,R)$ is normalized in $\Aut_R (R[X])$ by $\Gamma$. 
We thus have the \emph{$R$-general linear enrichment} of the action 
$\alpha$, given by 
$\mathcal{GL}(\Gamma,X,R) = \GL_f (X,R) \rtimes_{\hat{\alpha}} \Gamma$. 

The next result is due to Mal'cev \cite{Malc}; 
see Nica \cite{Nica} for a modern treatment. 

\begin{Thm}
If $R$ is a finitely generated integral domain then for every $n \geq 1$, 
$\GL_n(R)$ is residually finite. 
\end{Thm}

\begin{Cor} \label{LinearLocRFCoroll}
If $R$ is an integral domain, then $\GL_f (X,R)$ is locally residually finite, 
hence LEF, hence sofic. 
\end{Cor}

We deduce the following two results from Corollaries \ref{BigHaloCorollSpec} and \ref{LEFHaloProdCor}, respectively, 
with the final statement in each case following from Corollary \ref{LinearLocRFCoroll}. 

\begin{Thm} \label{LinEnrSoficThm}
    Suppose that $\Gamma$ is a sofic 
(respectively $\mathbb{K}$-linear sofic; hyperlinear; weakly sofic) 
group and that 
$\Gamma \curvearrowright_{\alpha} X$ is a sofic action on the set $X$. 
If $R$ is a unital ring with the property that 
$\Aut_R (R[Z])$ is a sofic ($\mathbb{K}$-linear sofic; hyperlinear; weakly sofic) group 
for every finite set $Z$, 
then $\mathcal{GL}(\Gamma,X,R)$ 
is a sofic ($\mathbb{K}$-linear sofic; hyperlinear; weakly sofic) group. 
In particular this holds if $R$ is an integral domain. 
\end{Thm}

\begin{Thm} \label{LinEnrLEFThm} 
    Suppose that $\Gamma$ is a LEF group and $\Gamma \curvearrowright_{\alpha} X$ is a LEF action on the set $X$. 
    If $R$ is a unital ring with the property that 
$\Aut_R (R[Z])$ is a LEF group for every finite set $Z$, 
then $\mathcal{GL}(\Gamma,X,R)$ is a LEF group. 
In particular this holds if $R$ is an integral domain. 
\end{Thm}

Concerning residual finiteness, 
we note that any nontrivial finitely-supported permutation of the set $X$ 
induces a nontrivial element of $\GL_f (X,R)$ by unique extension 
from a free basis to an $R$-module endomorphism of $R[X]$. 
In particular, if $X$ is infinite then $\GL_f (X,R)$ 
contains the infinite simple group $\Alt_f (X)$ as a subgroup. 
Meanwhile, if $R$ is an infinite field and $\lvert X \rvert \geq 2$, 
then $\GL_f (X,R)$ contains $\SL_2 (R)$ as a subgroup. 
In this setting it is well-known that every proper normal subgroup of $\SL_2 (R)$ 
is central, and $\PSL_2 (R)$ is an infinite simple group. 
Thus any non-central element of $\SL_2 (R)$ 
(viewed as an element of $\GL_f (X,R)$) survives in every finite quotient of $\GL_f (X,R)$. 
We summarise this discussion as follows. 

\begin{Prop} \label{ElemNoRFProp}
    Suppose $\mathcal{GL}(\Gamma,X,R)$ is a residually finite group. 
    Then $\Gamma$ is a residually finite group; $X$ is a finite set, 
    and if $\lvert X \rvert \geq 2$, $R$ is not an infinite field. 
\end{Prop}

In closing, we highlight two subgroups of $\mathcal{GL}(\Gamma,X,R)$ 
which have already received attention in the literature. 

For $x,y \in X$ distinct elements of the basis and $r \in R$, 
we have the \emph{elementary transvection} $E_{x,y} (r) \in \Aut_R (R[X])$, 
given by: 
\begin{center}
$E_{x,y} (r): x\mapsto x+ry$ and 
$E_{v,w} (r): z\mapsto z$ (for $x \neq z \in X$). 
\end{center}
Let $E_X (R)$ be the \emph{elementary subgroup} of $\Aut_R (R[X])$ 
generated by all $E_{x,y} (r)$. 
It is a classical fact that when $X$ is finite, 
$E_X (R) = \SL_{\lvert X \rvert} (R)$ for many rings $R$ 
(for instance this is the case when $R$ is a field or $R = \mathbb{Z}$). 

For a given ring $R$, the assignment $X \mapsto E_X (R)$ gives rise to 
a set-theoretic halo of groups, 
where for $Y \hookrightarrow X$ an inclusion of sets 
we can embed $E_Y (R)$ into $E_X (R)$ be declaring elements of 
$E_Y (R)$ to act trivially on the basis vectors $X \setminus Y$. 
Recalling that a group action $\Gamma \curvearrowright_{\alpha} X$ on the set $X$ 
induces an action $\hat{\alpha}$ by automorphisms on $\Aut_R (R[X])$, 
we note that under the latter action, $g$ maps $E_{x,y}(r)$ to 
$E_{\alpha(g)(x),\alpha(g)(y)}(r)$, so $E_X (R)$ is normalized 
in $\Aut_R (R[X])$ by $\Gamma$. 
We thus have the \emph{elementary enrichment} 
$\mathcal{E}(\Gamma,X,R) = E_X (R) \rtimes_{\hat{\alpha}} \Gamma$, 
previously studied in \cite{BradLEFExtn} and \cite{MimSak}. 

Genevois and Tessera \cite{GeneTess} consider instead the subgroup 
$\overline{E}_X (R)$ of $\Aut_R (R[X])$ generated by all $E_{x,y} (r)$ 
and all elements $D_x (\lambda) \in \Aut_R (R[X])$ given, for $x \in X$ 
and $\lambda \in R^{\ast}$, by: 
\begin{center}
$D_x (\lambda): x\mapsto \lambda x$ and 
$D_x (\lambda): z\mapsto z$ (for $x \neq z \in X$). 
\end{center}
Once again $X \mapsto \overline{E}_X (R)$ gives rise to a set-theoretic halo of groups, 
and $\overline{E}_X (R)$ is preserved under the action $\hat{\alpha}$ 
of $\Gamma$ on $\Aut_R (R[X])$, so we have an associated semidirect product 
$\overline{\mathcal{E}}(\Gamma,X,R)$, referred in \cite{GeneTess} as the 
\emph{lampcloner} group. 
Since $\mathcal{E}(\Gamma,X,R)$ and $\overline{\mathcal{E}}(\Gamma,X,R)$ embed 
as subgroups of $\mathcal{GL}(\Gamma,X,R)$, 
we have from Theorems \ref{LinEnrSoficThm} and \ref{LinEnrLEFThm} 
sufficient conditions for $\mathcal{E}(\Gamma,X,R)$ and $\overline{\mathcal{E}}(\Gamma,X,R)$ 
to be sofic; $\mathbb{K}$-linear sofic; hyperlinear; weakly sofic or LEF. 
In the case of LEF for $\mathcal{E}(\Gamma,X,R)$, 
we recover the qualitative part of the statement of Theorem 4.5 of \cite{BradLEFExtn} 
(see also Remark 4.15 of that paper). 

\subsection{Graph wreath products} \label{GraphWPSubsect}

Let $X=(V,E)$ be a graph and for each $v \in V$ let $H_v$ be a group. 
Write $\underline{H}$ to denote the tuple $(H_v)_{v \in V}$ 
of the groups $H_v$. 
Let $\ast_{v \in V} H_v$ be the free product of the $H_v$, and let: 
$$P(\underline{H},X) = (\ast_{v \in V} H_v)/\ll h^{-1} k^{-1} hk : (v,w) \in E , h \in H_v , k \in H_w \gg$$
be the \emph{graph-product} of the $H_v$. 
In the special case for which all the $H_v$ are isomorphic to a single group $H$, 
we shall write $P(H,X)$ for $P(\underline{H},X)$. 

Let $Y = (V',E')$ be an induced subgraph of $X$ 
and let $\iota : Y \hookrightarrow X$ be the inclusion map. 
Using the universal property of free products, 
we have inclusion and projection homomorphism: 
$$\tilde{\iota} :\ast_{v \in V'} H_v \rightarrow\ast_{v \in V} H_v 
\text{ and } \tilde{\pi} :\ast_{v \in V} H_v \rightarrow\ast_{v \in V'} H_v $$
where $\tilde{\pi}$ sends $H_v$ identically to itself for each $v \in V'$, 
and $\tilde{\pi}(H_v) = \lbrace e \rbrace$ for all $v \in V \setminus V'$. 
Then $\tilde{\pi} \circ \tilde{\iota}$ 
is the identity map on $\ast_{v \in V'} H_v$, 
and since $Y$ is an induced subgraph of $X$, 
$\tilde{\iota}$ and $\tilde{\pi}$ descend to well-defined homomorphisms 
$\hat{\iota} : P(\underline{H},Y) \rightarrow P(\underline{H},X)$ 
and $\hat{\pi} : P(\underline{H},X) \rightarrow P(\underline{H},Y)$ satisfying
$\hat{\pi} \circ \hat{\iota} = \id_{P(\underline{H},Y)}$, 
so that $\hat{\iota}$ is a monomorphism. 
We conclude the following. 

\begin{Lemma} \label{GraphWPEmbedLemma}
    For any induced subgraph $Y = (V',E')$ of $X$, 
    the subgroup of $P(\underline{H},X)$ generated by the image of $\bigcup_{v \in V'} H_v$ 
    is isomorphic to $P(\underline{H},Y)$. 
    In particular, the image of $H_v$ in $P(\underline{H},X)$ 
    is isomorphic to $H_v$, and for $v , w \in V$, 
    the subgroup of $P(\underline{H},X)$ generated by the image of $H_v \cup H_w$ 
    is isomorphic to $H_v \times H_w$ (respectively $H_v \ast H_w$) 
    if $(v,w) \in E$ (respectively $(v,w) \notin E$). 
\end{Lemma}

Henceforth, for $Y$ an induced subgraph of $X$, 
we shall identify $P(\underline{H},Y)$ with its image in $P(\underline{H},X)$ under 
$\hat{\iota}$ without further comment. 

We also recall a few notions from Chapter 3 of \cite{GreenThesis} 
to be used in the remainder of the Subsection. 
By a \emph{word} $w$ in $P(\underline{H},X)$ we shall mean an 
expression $w = h_1 \cdots h_n$, with $e \neq h_i \in H_{v_i}$ 
for some vertices $v_1 , \ldots , v_n \in V$. 
The $h_i$ are the \emph{syllables} of the word $w$. 
If $v_i$ and $v_{i+1}$ are distinct but adjacent vertices in $X$, then:
$$h_1 \cdots h_i h_{i+1} \cdots h_n = h_1 \cdots h_{i+1} h_i \cdots h_n$$
in $P(\underline{H},X)$. A finite sequence of such transformations 
will be referred to as a \emph{syllable-shuffle} of the word $w$. 
A \emph{reduction} of a word $w$ consists of syllable-shuffling $w$ such that 
for some $i$ we have $v_i = v_{i+1}$, then either replacing the two syllables 
$h_i$, $h_{i+1}$ by the single syllable $h_i h_{i+1}$ (if $h_{i+1} \neq h_i ^{-1}$) 
or deleting both (if $h_{i+1} = h_i ^{-1}$). 
A word is \emph{reduced} if no reductions can be performed on it. 
Every element of $P(\underline{H},X)$ is represented by a reduced word, 
which is unique up to syllable-shuffling (\cite{GreenThesis} Theorem 3.9). 
Since applying a reduction to $w = h_1 \cdots h_n$ reduces the value of $n$, 
any sequence of reductions applied to a word will eventually lead to a reduced word. 
If $w = h_1 \cdots h_n$ is reduced, then $\lbrace v_1 , \ldots , v_n \rbrace$ 
is the \emph{support} of $w$, and $n$ is the \emph{syllable-length} of $w$. 

\begin{Prop} \label{GWPHaloProp}
For any fixed group $H$ the assignment $X \mapsto P(H,X)$ 
induces a graph-theoretic halo of groups. 
\end{Prop}

\begin{proof}
We have already seen above that an inclusion of an induced subgraph induces 
an embedding of graph-products. The functoriality of this embedding is easy to verify. 
Axioms (i) and (ii) from Definition \ref{HaloDefn} are easy to verify; 
for (iii) suppose that $Y$ and $Z$ are induced subgraphs of $X$, 
and there are $k_Y \in P(H,Y)$, $k_Z \in P(H,Z)$ such that 
$k_Y = k_Z$ in $P(H,X)$. 
We may write $k_Y$ as a reduced word with support in $V(Y)$ and 
$k_Z $ as a reduced word with support in $V(Z)$, but by uniqueness of reduced words 
(up to syllable-shuffling) these reduced words have the same support, 
which is therefore in $V(Y \cap Z)$. 
\end{proof}

If $\Gamma \curvearrowright_{\alpha} X$ is a group action on $X$ 
by graph automorphisms, then there is an induced action $\hat{\alpha}$ of 
$\Gamma$ on $P(H,X)$ by group automorphisms, 
as per Definition \ref{FunctHaloProdDefn}. 
The associated semidirect product is the \emph{graph wreath product} 
$\mathcal{W}(\Gamma,H,X) = P(H,X) \rtimes_{\hat{\alpha}} \Gamma$. 
We highlight two extreme cases which have received particular attention. 
\begin{itemize}
\item[(i)] If $X$ is a complete graph, then $P(H,X)$ 
is the direct sum of copies of $H$ indexed by $V$, 
and $\mathcal{W}(\Gamma,H,X)$ is the (restricted permutational) 
\emph{wreath product} $H\wr_{\alpha} \Gamma$ of $H$ and $\Gamma$; 

\item[(ii)] If $X$ is the empty graph on $V$, 
then $P(H,X)$ is the free product of copies of $H$ indexed by $V$, 
and $\mathcal{W}(\Gamma,H,X)$ is the \emph{free wreath product} 
$H\wr_{\alpha} ^{\ast} \Gamma$. 
In the special case for which $\alpha$ is the regular action 
of $\Gamma$ on itself, the free wreath product of $H$ and $\Gamma$ is 
easily seen to be isomorphic to the free product $H \ast \Gamma$ 
of $H$ and $\Gamma$. 
\end{itemize}

Of course, since every induced subgraph of a complete (respectively empty) 
graph is complete (empty), these special cases give rise to 
\emph{set-theoretic} haloes of groups. 

Metric approximations of wreath products was studied by Hayes and Sale 
\cite{HayesSaleSofic,HayesSaleMetric}, 
who proved that the regular restricted wreath product 
(that is, the case of $H\wr_\alpha \Gamma$ for which 
$\alpha$ is the regular action of $\Gamma$ on itself) 
of two sofic (respectively $\mathbb{K}$-linear sofic; hyperlinear; weakly sofic) groups is 
sofic ($\mathbb{K}$-linear sofic; hyperlinear; weakly sofic). 
Elek and Szabo \cite{ElekSzaboSofic} proved that 
the class of sofic groups is closed under free products. 
These results were greatly generalized by Gao, Elayavalli and Patchell \cite{GaoElayPatchSet,GaoElayPatchGraph}, 
to all graph wreath products and to all sofic actions. 

\begin{Thm} \label{GraphWPSoficThm}
Suppose that $\Gamma$ and $H$ are sofic groups, and that 
$\Gamma \curvearrowright_{\alpha} X$ is a sofic action on the graph $X$. 
Then $\mathcal{W}(\Gamma,H,X)$ is a sofic group. 
\end{Thm}

Note that, by the above discussion, Theorem \ref{GraphWPSoficThm} 
also follows from our Corollary \ref{BigSemiCoroll} and the 
fact that for any sofic group $H$ and any graph $X$, 
$P(H,X)$ is a sofic group (a result of Ciobanu, Holt and Rees \cite{CiobHoltRees}). 

Gao, Elayavalli and Patchell \cite{GaoElayPatchSet,GaoElayPatchGraph} 
provide a similar result for hyperlinearity. 

\begin{Thm} \label{GWPHypThm}
Suppose that $\Gamma$ and $H$ are hyperlinear groups, and that 
$\Gamma \curvearrowright_{\alpha} X$ is a sofic action on the graph $X$. 
Then $\mathcal{W}(\Gamma,H,X)$ is a hyperlinear group. 
\end{Thm}

As in the sofic case, we may recover their result, 
in this case using the following result of Caspers \cite{Caspers}. 

\begin{Thm} \label{CaspersThm}
The graph product of hyperlinear groups is hyperlinear.
\end{Thm}

\begin{proof}
A group is hyperlinear iff its group von Neumann algebra
is Connes embeddable. 
The result follows from Theorem 3.4 of \cite{Caspers},
since the von Neumann algebra of a graph product of
groups $\Gamma_v$ is isomorphic to the
corresponding graph product of the von Neumann algebras
of the $\Gamma_v$ (as noted in the Introduction to \cite{Caspers}). 
\end{proof}

\begin{proof}[Proof of Theorem \ref{GWPHypThm}]
This is immediate from Theorem \ref{CaspersThm}; Proposition \ref{GWPHaloProp} 
and Corollary \ref{BigSemiCoroll}. 
\end{proof}

Once again, Theorem \ref{GWPHypThm} was originally obtained as Theorem 1.12 of \cite{GaoElayPatchGraph} 
(again recalling that a group is hyperlinear iff its group von Neumann algebra
is Connes embeddable). 

For linear soficity, the state of knowledge is less well-developed. 
Stolz (Theorem 5.6 of \cite{Stolz}) has claimed that for a field $\mathbb{K}$, 
the class of $\mathbb{K}$-linear sofic groups is closed under free products. 
Alas, there is a gap in the proof of this Theorem, as detailed in \cite{NeedhamErratum}. 
It is easy to see that the direct sum of the class of $\mathbb{K}$-linear sofic groups is 
closed under direct sums, but at present this seems to be all that is known 
about the permanence of $\mathbb{K}$-linear soficity under graph products. 

\begin{Conj} \label{GWPLinSoficConj}
Let $\mathbb{K}$ be a field; let $V$ be a set and let $\lbrace \Gamma_v : v \in V \rbrace$ 
be a family of $\mathbb{K}$-linear sofic groups. Then: 
\begin{itemize}
    \item[(a)] $\ast_{v \in V} \Gamma_v$ is $\mathbb{K}$-linear sofic; 

    \item[(b)] More generally, for any graph $X = (V,E)$ on the vertex set $V$, 
    $$(\ast_{v \in V} \Gamma_v)/\ll h^{-1} k^{-1} hk : (v,w) \in E , h \in \Gamma_v , k \in \Gamma_w \gg$$
    is $\mathbb{K}$-linear sofic. 
\end{itemize}
\end{Conj}

\begin{Thm}
Let $\mathbb{K}$ be a field. Suppose that $\Gamma$ and $H$ are $\mathbb{K}$-linear sofic groups; 
that $\Gamma \curvearrowright_{\alpha} X$ is a sofic action on the set $X$, 
and that $\Gamma \curvearrowright_{\beta} Y$ is a sofic action on the graph $y$. Then: 
\begin{itemize}
    \item[(a)] $H \wr_{\alpha} \Gamma$ is $\mathbb{K}$-linear sofic; 
    
    \item[(b)] If Conjecture \ref{GWPLinSoficConj} (a) is true, 
    then $H \wr_{\alpha} ^{\ast} \Gamma$ is $\mathbb{K}$-linear sofic; 
    
    \item[(c)] If Conjecture \ref{GWPLinSoficConj} (b) is true, then 
    $\mathcal{W}(\Gamma,H,Y)$ is $\mathbb{K}$-linear sofic. 
\end{itemize}
\end{Thm}

\begin{proof}
This is immediate from Proposition \ref{GWPHaloProp} 
and Corollary \ref{BigSemiCoroll}. 
\end{proof}

For weakly sofic groups, the situation is very similar to the 
$\mathbb{K}$-linear sofic case. A direct sum of weakly sofic groups 
is easily seen to be weakly sofic, but the corresponding problem for other graph products 
(and even free products) appears to be open. 

\begin{Conj} \label{GWPWSoficConj}
Let $V$ be a set and let $\lbrace \Gamma_v : v \in V \rbrace$ 
be a family of weakly sofic groups. Then: 
\begin{itemize}
    \item[(a)] $\ast_{v \in V} \Gamma_v$ is weakly sofic; 

    \item[(b)] More generally, for any graph $X = (V,E)$ on the vertex set $V$, 
    $$(\ast_{v \in V} \Gamma_v)/\ll h^{-1} k^{-1} hk : (v,w) \in E , h \in \Gamma_v , k \in \Gamma_w \gg$$
    is weakly sofic. 
\end{itemize}
\end{Conj}

\begin{Thm} \label{GraphWPWSoficThm}
Suppose that $\Gamma$ and $H$ are weakly sofic groups; 
that $\Gamma \curvearrowright_{\alpha} X$ is a sofic action on the set $X$, 
and that $\Gamma \curvearrowright_{\beta} Y$ is a sofic action on the graph $y$. Then: 
\begin{itemize}
    \item[(a)] $H \wr_{\alpha} \Gamma$ is weakly sofic; 
    
    \item[(b)] If Conjecture \ref{GWPWSoficConj} (a) is true, 
    then $H \wr_{\alpha} ^{\ast} \Gamma$ is weakly sofic; 
    
    \item[(c)] If Conjecture \ref{GWPWSoficConj} (b) is true, then 
    $\mathcal{W}(\Gamma,H,Y)$ is weakly sofic. 
\end{itemize}
\end{Thm}

\begin{proof}
This is immediate from Proposition \ref{GWPHaloProp} 
and Corollary \ref{BigSemiCoroll}. 
\end{proof}

Regarding LEF for graph wreath products, 
Vershik and Gordon \cite{VersGord} proved that the regular restricted 
wreath product of two LEF groups is LEF. 
It is rather easy to see that the free product of 
two LEF groups is LEF, using the fact that a free product 
of finite groups is residually finite. 
More generally, we have the following. 

\begin{Thm} \label{GraphProdLEFThm}
For any graph $X=(V,E)$ and any LEF groups $(H_v)_{v \in V}$, $P(\underline{H},X)$ is an LEF group. 
\end{Thm}

\begin{proof}
Let $F \subseteq P(\underline{H},X)$ be finite. 
First, there is a finite induced subgraph $Y$ of $X$ 
such that $F \subseteq P(\underline{H},Y)$. 
Thus the projection homomorphism $\hat{\pi} : P(\underline{H},X) \rightarrow P(\underline{H},Y)$ 
defined above restricts to an injection on $F$. 
Replacing $X$ by $Y$, we have reduced to the case for which $X$ is a finite graph. 

Second, there exists $L \in \mathbb{N}$ and finite subsets $F_v \subseteq H_v$ 
such that every element of $F$ can be written as a reduced word $w$ 
of syllable-length at most $L$, 
such that for each $v \in V$, all the $H_v$-syllables of $w$ lie in $F_v$. 
Enlarging $F_v$ if required, we may assume $e \in F_v = F_v ^{-1}$ for all $v \in V$. 
Let $F_v ' \subseteq H_v$ consist of all products of at most $2 L$ elements of $F_v$. 
Since $H_v$ is a LEF group, there exists a finite group $Q_v$ and a local embedding 
$\phi_v : F_v ' \rightarrow Q_v$. 
Define $\Phi : F \rightarrow P(\underline{Q},X)$ such that for $w = h_1 \cdots h_n$ 
a reduced word with $e \neq h_i \in H_{v_i}$, 
$\Phi(w) = \phi_{v_1} (h_1) \cdots \phi_{v_n} (h_n)$. 
Note that, since $P(\underline{H},X)$ and $P(\underline{Q},X)$ 
are defined on the same graph $X$, $\phi_{v_1} (h_1) \cdots \phi_{v_n} (h_n)$ 
is also a reduced word, and that a syllable-shuffling of $w$ 
induces a syllable-shuffling of $\Phi(w)$, and vice versa. 
Thus $\Phi$ is a well-defined injective function. 

We check that $\Phi$ is a local embedding. 
Let $w = h_1 \cdots h_n$ and $\overline{w} = \overline{h}_1 \cdots \overline{h}_m$ 
be reduced words in $F$ (with $v_i , \overline{v}_j \in V$ 
such that $e \neq h_i \in H_{v_i}$ and $e \neq \overline{h}_i \in H_{\overline{v}_i}$) 
and suppose $w \overline{w} \in F$. 
Let $h_1 ' \cdots h_l '$ be a reduced word representing $w \overline{w}$ 
(with $v_i ' \in V$ such that $e \neq h_i '\in H_{v_i '}$), so that: 
$$\Phi(w) \Phi(\overline{w}) = \phi_{v_1} (h_1) \cdots \phi_{v_n} (h_n)\phi_{\overline{v}_1} (\overline{h}_1) \cdots \phi_{\overline{v}_m} (\overline{h}_m)$$
and 
$$\Phi(w \overline{w}) = \phi_{v_1 '} (h_1 ') \cdots \phi_{v_n '} (h_n ')$$
so that it suffices to check that $\phi_{v_1 '} (h_1 ') \cdots \phi_{v_n '} (h_n ')$ 
is a reduced word representing $\Phi(w) \Phi(\overline{w})$. 
To this end, consider a sequence of reductions which transforms 
$h_1 \cdots h_n \overline{h}_1 \cdots \overline{h}_m$ into 
$h_1 ' \cdots h_l '$. 
After $i$ reductions, we have a word $h_1 ^{(i)} \cdots h_{l_i} ^{(i)}$ 
(with $v_j ^{(i)} \in V$ such that $e \neq h_j ^{(i)} \in H_{v_j ^{(i)}}$) 
with $l_i \leq n+m-i$. 
Suppose by induction that: 
$$\Phi(w) \Phi(\overline{w}) = \phi_{v_1 ^{(i)}} (h_1 ^{(i)}) \cdots \phi_{v_{l_i} ^{(i)}} (h_{l_i} ^{(i)})$$
and that each $h_j ^{(i)}$ 
is a product of at most $n_j$ elements of $F_{v_j ^{(i)}}$, such that 
$n_1 + \ldots + n_{l_i} \leq 2 L$. 
Up to syllable-shuffling, we may assume that for some 
$1 \leq j \leq l_i - 1$ we have $v_j ^{(i)} = v_{j+1} ^{(i)} = v^{(i)}$. 
Since $h_j ^{(i)}$, $h_{j+1} ^{(i)}$ and $h_j ^{(i)} h_{j+1} ^{(i)}$ 
are elements of $F_{v^{(i)}} '$, we have:  
\begin{equation} \label{LEFGraphProdEqn}
\phi_{v^{(i)}} (h_j ^{(i)} h_{j+1} ^{(i)})=\phi_{v^{(i)}} (h_j ^{(i)})\phi_{v^{(i)}} (h_{j+1} ^{(i)}).
\end{equation}
If $(h_j ^{(i)}) ^{-1} = h_{j+1} ^{(i)}$ then 
we apply a reduction to delete $h_j ^{(i)} $ and $h_{j+1} ^{(i)}$ 
to obtain a new word $h_1 ^{(i)} \cdots h_{j-1} ^{(i)} h_{j+2} ^{(i)} \cdots h_{l_i} ^{(i)}$ 
representing $w \overline{w}$, 
and by (\ref{LEFGraphProdEqn}), 
$$\Phi(w) \Phi(\overline{w}) 
= \phi_{v_1 ^{(i)}} (h_1 ^{(i)}) \cdots  \phi_{v_{j-1} ^{(i)}} (h_{j-1} ^{(i)})\phi_{v_{j+2} ^{(i)}} (h_{j+2} ^{(i)}) \cdots \phi_{v_{l_i} ^{(i)}} (h_{l_i} ^{(i)})$$
so that we may continue our inductive argument. 
Similarly if $(h_j ^{(i)}) ^{-1} \neq h_{j+1} ^{(i)}$ then our reduction 
combines $h_j ^{(i)}$ and $h_{j+1} ^{(i)}$ into $h_j ^{(i+1)} = h_j ^{(i)} h_{j+1} ^{(i)}$ 
(a product of at most $n_j + n_{j+1}$ elements of $F_{v^{(i)}}$); 
our new word representing $w \overline{w}$ is: 
$$h_1 ^{(i)} \cdots h_{j-1} ^{(i)} h_j ^{(i+1)} h_{j+2} ^{(i)} \cdots h_{l_i} ^{(i)}$$
and by (\ref{LEFGraphProdEqn}) we have
$$\Phi(w) \Phi(\overline{w}) 
= \phi_{v_1 ^{(i)}} (h_1 ^{(i)}) \cdots  \phi_{v_{j-1} ^{(i)}} (h_{j-1} ^{(i)}) 
\phi_{v^{(i)}} (h_j ^{(i+1)} )
\phi_{v_{j+2} ^{(i)}} (h_{j+2} ^{(i)}) \cdots \phi_{v_{l_i} ^{(i)}} (h_{l_i} ^{(i)})$$
so that our induction may continue in this case also. 
Thus $\Phi$ is indeed a local embedding. 

Finally, by a result of Green (Corollary 5.4 of \cite{GreenThesis}), 
$P(\underline{Q},X)$ is residually finite, so there is a finite group $\overline{Q}$ 
and a homomorphism $\rho : P(\underline{Q},X) \rightarrow \overline{Q}$ 
with $\rho |_{\Phi(F)}$ injective. 
Then $\rho \circ \Phi$ is the desired injective partial homomorphism. 
\end{proof}

We deduce: 

\begin{Thm} \label{GraphWPLEFThm}
For $\Gamma$, $H$ groups, with $H$ nontrivial, and $X$ a graph, 
the graph wreath product $\mathcal{W}(\Gamma,H,X)$ is a LEF group 
iff $\Gamma$ and $H$ are LEF groups and 
$\Gamma \curvearrowright_{\alpha} X$ is an LEF action on the graph $X$. 
\end{Thm}

\begin{proof}
For the reverse direction, 
by Proposition \ref{GWPHaloProp}, Corollary \ref{LEFHaloProdCor} is applicable. 
The desired conclusion is then immediate from Theorem \ref{GraphProdLEFThm}. 

For the forward direction, it is clearly necessary that 
$\Gamma$ and $H$ be LEF groups, since both are isomorphic 
to subgroups of $\mathcal{W}(\Gamma,H,X)$. 
Let $X=(V,E)$ and let $Z \subseteq V$, $F \subseteq \Gamma$ be finite. 
Set $e \neq h \in \Delta$ and for $v \in V$ 
let $h_v\in H_v \leq\mathcal{W}(\Gamma,H,X)$ 
be the copy of $h$ supported at vertex $v$. 
Let $F' \subseteq \mathcal{W}(\Gamma,H,X)$ 
be a finite subset containing as subsets: $e \cup F \cup F^{-1}$, 
$\lbrace h_z , h_z h_y , h_y h_z : y,z \in Z \rbrace$ 
and $\lbrace g h_z , h_z g^{-1} : z \in Z , g \in F \rbrace$. 
Let $\phi : F' \hookrightarrow Q$ be an injective partial 
homomorphism, for some finite group $Q$. 

Define a graph $\overline{X} = (\overline{V},\overline{E})$ 
with $\overline{V} = Q$, with two vertices being adjacent 
iff the corresponding elements of $Q$ commute. 
Then the action of $Q$ by conjugation on itself clearly 
induces an action of $Q$ on $\overline{X}$ by graph automorphisms. 
Define $\psi : F \rightarrow Q$ by restriction of $\phi$ 
and let $\pi : Z \rightarrow \overline{V}$ be given by 
$\pi(z) = \phi (h_z)$. 
Then $\psi$ is a partial homomorphism; 
$\pi$ is the inclusion of an induced subgraph 
(since by Lemma \ref{GraphWPEmbedLemma}, $h_v$ and $h_w$ commute in $\mathcal{W}(\Gamma,H,X)$ 
iff $v,w \in V$ are adjacent or equal), 
and for $g\in F$, $z\in Z$ such that $\alpha(g)\cdot z\in Z$, 
we have:
$$\pi (\alpha(g) \cdot z) = \phi (h_{\alpha(g) \cdot z}) 
 = \phi (g h_z g^{-1}) = \phi (g) \phi (h_z) \phi (g)^{-1}
 = \psi (g) \pi (z) \psi (g)^{-1}$$
so that the conditions of Definition \ref{LEFActEAAutoProp} are satisfied. 
\end{proof}

The special case of our Theorem \ref{GraphWPLEFThm} for permutational 
wreath products is due to Cornulier \cite{CornulierWreath}. 

\begin{Prop}[\cite{CornulierWreath} Proposition 1.5] \label{CornWreathProp}
Let $\Gamma$ and $H$ be groups, and let $\Delta \leq \Gamma$. 
Then $H \wr_{\Gamma/\Delta} \Gamma$ is LEF iff 
$\Gamma$ and $H$ are LEF and $\Delta$ is AFI in $\Gamma$. 
\end{Prop}


Residual finiteness of graph wreath products has been characterised by Needham \cite{Needham}. 
For $X=(V,E)$ a graph and $v \in V$ let $N(v) \subseteq V \setminus \lbrace v \rbrace$ denote the set of neighbours or $v$ in $X$. 

\begin{Thm}[Needham] \label{NeedhamGWPRFThm}
For $\Gamma$ and $H$ groups; $X=(V,E)$ a graph and 
$\Gamma \curvearrowright_{\alpha} X$ an action by graph automorphisms, 
graph wreath product $\mathcal{W}(\Gamma,H,X)$ 
is residually finite iff the following hold. 
\begin{itemize}
    \item[(i)] $\Gamma$ and $H$ are residually finite; 
    \item[(ii)] Either $H$ is abelian and for all $(v,w) \in E$ 
    there exists $K \leq _f \Gamma$ such that $w \notin \alpha(K)(v)$, 
    or for all $v \in V$ there exists $K \leq _f \Gamma$ such that $\alpha(K)(v)$ 
    contains no neighbour of $v$; 
    \item[(iii)] For all distinct vertices $v,w \in V$ such that $v,w \notin E$, 
    there exists $K \leq _f \Gamma$ such that $\alpha(K)(w)$ does not contain $v$ 
    or any neighbour of $v$. 
\end{itemize}
\end{Thm}

\subsection{Verbal wreath products}

We follow the treatment of Brude and Sasyk \cite{BrudSasy}.
Let $F_{\infty}$ be a free group of countably infinite rank
on a basis $x_1 , x_2 , \ldots$,
and let $W \subseteq F_{\infty}$ be a set of words.
For $X$ be a set and $H$ be a fixed group,
let $F(H,X) = \ast_{x \in X} H_x$ be the free product
of groups $H_x \cong H$, indexed by the elements of $X$.
Let $W(F(H,X)) \vartriangleleft F(H,X)$
be the verbal subgroup corresponding to the set of words $W$,
and let $[H]^{F(H,X)}$ be the normal subgroup of $F(H,X)$
generated by all elements of the form
$h_x ^{-1} k_y ^{-1} h_x k_y$ for $x,y \in X$
with $x \neq y$ and $h_x \in H_x$, $k_y \in H_y$.
Then the \emph{verbal product} of the $H_x$ is defined to be:
$$F^W (H,X) = F(H,X) / \big(W(F(H,X)) \cap [H]^{F(H,X)}\big)$$

If $\iota : Y \hookrightarrow X$ is an injection of sets,
then (by the universal property of free products),
there is an induced monomorphism
$\tilde{\iota} : F(H,Y) \hookrightarrow F(H,X)$,
sending each $H_y$ to $H_{\iota (y)}$.
It is clear that:
$$\tilde{\iota} \big( W(F(H,Y)) \big) \leq W(F(H,X))
\text{ and } \tilde{\iota} \big( [H]^{F(H,Y)} \big) \leq [H]^{F(H,X)}$$
so that $\tilde{\iota}$ descends to a well-defined
homomorphism $\hat{\iota} : F^W (H,Y) \rightarrow F^W (H,X)$.

\begin{Prop} \label{VerbProdHaloProp}
For each group $H$ and each $W \subseteq F_{\infty}$,
the assignment $X \mapsto F^W (H,X)$ induces a set-theoretic
halo of groups.
\end{Prop}

\begin{proof}
It was already noted in Section 2 of \cite{GeneTess} that,
owing to Lemma 4.6 of \cite{BrudSasy},
$F^W (H,X)$ gives rise to a halo of groups.
From Lemma 4.4 of \cite{BrudSasy}, we have that the
homomorphism $\hat{\iota} : F^W (H,Y) \rightarrow F^W (H,X)$
induced by an injection $\iota : Y \hookrightarrow X$
is a monomorphism (we may choose an arbitrary total ordering on $X$).
The functoriality properties of $\hat{\iota}$
are clear from the definition of $\hat{\iota}$ given above.
\end{proof}

For a group action $\Gamma \curvearrowright_{\alpha} X$ on a set $X$,
we therefore have an induced action
$\Gamma \curvearrowright_{\hat{\alpha}} F^W (H,X)$
by automorphisms
hence we may form the associated \emph{verbal wreath product}
$H \wr^W _{\alpha}  \Gamma = F^W (H,X) \rtimes_{\hat{\alpha}} \Gamma$.

Brude and Sasyk give particular attention to the following sets
of words $W \subseteq F_{\infty}$, and to the associated
verbal product groups.

\begin{Ex} \label{VerbalProdExDefn}
Let $H$ be a group.
\begin{itemize}
\item[(i)] If $W = \lbrace n_k \rbrace$,
where $n_k \in F_{\infty}$ ($k \in \mathbb{N}$)
is defined inductively by
$n_1 = [x_2,x_1]$ and $n_k = [x_{k+1},n_k]$,
then $F^W (H,X)$ is the \emph{$k$-nilpotent product} of the $H_x$;

\item[(ii)] If $W = \lbrace s_k \rbrace$,
where $s_k \in F_{\infty}$ ($k \in \mathbb{N}$)
is defined inductively by
$s_1 (x_1,x_2) = [x_2,x_1]$ and:
$$s_k (x_1 , \ldots , x_{2^k})
= [s_{k-1}(x_1 , \ldots , x_{2^{k-1}}),s_{k-1}(x_{2^{k-1} + 1} , \ldots , x_{2^k})],$$
then $F^W (H,X)$ is the \emph{$k$-soluble product} of the $H_x$;

\item[(iii)] If $W = \lbrace x_1 ^k \rbrace$, for $k \in \mathbb{N}$,
then $F^W (H,X)$ is the \emph{$k$-Burnside product} of the $H_x$

\end{itemize}
\end{Ex}

The next result is the content of Theorems 1.3-1.5 of \cite{BrudSasy}.

\begin{Thm} \label{BrudSasyProdThm}
Let $\mathcal{C}$ be the class of sofic (respectively
$\mathbb{K}$-linear sofic; hyperlinear; weakly sofic) groups.
Let $W \subseteq F_{\infty}$ be one of the sets from
Example \ref{VerbalProdExDefn}
(with the stipulation in case (iii) that $k=2,3,4$ or $6$).
If $H \in \mathcal{C}$, then for any set $X$,
$F^W (H,X) \in \mathcal{C}$.
\end{Thm}

We deduce: 

\begin{Thm} \label{MainVerbalProdThm}
Let $\Gamma \curvearrowright_{\alpha} X$ be a sofic
action of the group $\Gamma$ on the set $X$.
Let $\mathcal{C}$ be the class of sofic (respectively
$\mathbb{K}$-linear sofic; hyperlinear; weakly sofic) groups.
Suppose that $H \in \mathcal{C}$ and that $W \subseteq F_{\infty}$
is such that $F^W (H,Z) \in \mathcal{C}$ for all finite sets $Z$.
If $\Gamma \in \mathcal{C}$
then $H \wr^W _{\alpha} \Gamma \in \mathcal{C}$.
In particular, if $W$ is as in Theorem \ref{BrudSasyProdThm},
then $H \wr^W _{\alpha}  \Gamma \in \mathcal{C}$ whenever $H , \Gamma \in \mathcal{C}$.
\end{Thm}

The special case of Theorem \ref{MainVerbalProdThm}
for which $\Gamma$ is a sofic group and $\alpha$
is the regular action of $\Gamma$ on itself is
the content of Theorems 1.8 and 1.9 of \cite{BrudSasy}.

\begin{proof}[Proof of Theorem \ref{MainVerbalProdThm}]
By Proposition \ref{VerbProdHaloProp} and Theorem \ref{BrudSasyProdThm},
Corollary \ref{BigHaloCorollSpec} is applicable.
\end{proof}

The situation regarding LEF and residual finiteness
for verbal wreath products seems much less clear. 
Golovin \cite{Golov} proved the following. 

\begin{Thm} \label{GolovinThm}
    For each $k \geq 1$ the $k$-nilpotent product of a finite family of finite groups is finite. 
\end{Thm}

We deduce: 

\begin{Thm} \label{LEFVerbalProdThm}
    Let $k \geq 1$ and let $W = \lbrace n_k \rbrace$ be as in Example \ref{VerbalProdExDefn} (i). 
Let $\Gamma \curvearrowright_{\alpha} X$ be a LEF 
action of the group $\Gamma$ on the set $X$, and let $H$ be a finite group. 
Then $H \wr^W _{\alpha}  \Gamma$ is a LEF group. 
\end{Thm}

\begin{proof}
By Theorem \ref{GolovinThm} and Proposition \ref{VerbProdHaloProp}, 
$F^W (H,X)$ is locally finite, hence LEF. 
The conclusion then follows from Corollary \ref{LEFHaloProdCor}. 
\end{proof}

The conclusion of Theorem \ref{LEFVerbalProdThm} could be extended to 
cover the case of $H$ a LEF group, conditional on the 
following very plausible extension of Golovin's Theorem. 

\begin{Conj}
    For each $k \geq 1$ the $k$-nilpotent product of a finite family of LEF groups is LEF. 
\end{Conj}

\subsection{Automorphic enrichments}

Let $\mathcal{L}_X = \lbrace L(Y) : Y \subseteq X \rbrace$ be a  halo of groups over the set $X$ and let: 
$$\AL_f (X) = \lbrace \beta \in \Aut(L(X)) : \exists Y \subseteq X \text{ finite s.t. } \beta(L(Y))=L(Y) \text{ and } \beta_{L(X \setminus Y)} = \id_{L(X \setminus Y)} \rbrace$$
and for $Y \subseteq X$ let: 
$$\AL_{f,X} (Y) = \lbrace \beta \in \AL_f (X) : \beta(L(Y))=L(Y) \text{ and } \beta_{L(X \setminus Y)} = \id_{L(X \setminus Y)} \rbrace$$
so that $\AL_f (X)$ is the directed union of the directed system 
$\mathcal{AL}_X = \lbrace \AL_{f,X} (Y) : Y \subseteq X \rbrace$ and is generated by 
$\bigcup_{Z \subseteq X \text{ finite}} \AL_{f,X} (Z)$. 

\begin{Rem} \label{AutMostlyHalo}
It is clear that $\mathcal{AL}_X$ satisfies items (i) and (ii) of Definition \ref{HaloDefn}, 
and that for $Y,Z \subseteq X$, every $\beta \in \AL_{f,X} (Y) \cap \AL_{f,X} (Z)$ 
preserves $L(Y) \cap L(Z) = L(Y \cap Z)$. 
\end{Rem}

If $\Gamma \curvearrowright_{\alpha} X$ is a group action 
and $\Gamma \curvearrowright_{\hat{\alpha}} L(X)$ 
is an action by automorphsms, such that: 
$$\hat{\alpha}(g) \big( L(Y) \big) = L(\alpha(g) (Y))$$
for all $g \in \Gamma$ and $Y \subseteq X$, 
we have an induced action of $\Gamma$ by automorphisms on $\AL_{X,f} (X)$, 
and may therefore define the 
\emph{$\mathcal{L}_X$-automorphic enrichment 
of $\Gamma \curvearrowright_{\alpha} X$} to be 
$\mathcal{AL}(\Gamma,X) = \AL_f (X) \rtimes_{\hat{\alpha}} \Gamma$. 

\begin{Ex}
In \cite{GeneTess} the following instance of the above construction 
was described. Let $F(X)$ be the free group with basis the set $X$. 
Then the family $\lbrace F(Y) : Y \subseteq X \rbrace$ 
of subgroups of $F(X)$ forms a halo of groups over $X$. 
It is clear that any group action $\Gamma \curvearrowright_{\alpha} X$ 
extends to an action by automorphisms on $F(X)$ such that 
for all $g \in \Gamma$ and $Y \subseteq X$, 
$\alpha(g) (F(Y)) = F(\alpha(g) (Y))$. 
We thus have an action of $\Gamma$ by automorphisms on 
$\FAut (F(X))$ (the group of automorphisms of $F(X)$ 
supported on some finite-rank subgroup of $F(X)$), 
and an associated halo product $\FAut (F(X)) \rtimes_{\hat{\alpha}} \Gamma$. 
\end{Ex}

It shall be useful to consider haloes of groups satisfying 
two additional structural properties. 

\begin{itemize}
\item[(a)] For all $Y,Z \subseteq X$, 
$L(Y \cup Z) = \langle L(Y) \cup L(Z) \rangle$; 

\item[(b)] If $Y \subseteq X$ with $Y$ finite, 
then for all $\beta \in \Aut(L(Y))$ 
there exists $\tilde{\beta} \in \AL_{f,X} (Y)$ extending $\beta$. 
\end{itemize}

\begin{Rem} \label{AutOnePtGenRem}
It is easy to see that property (a) is equivalent to: 
for all $Y \subseteq X$, 
$$L(Y) = \langle L(\lbrace y \rbrace) : y \in Y \rangle.$$
\end{Rem}

\begin{Ex} \label{SumProdHaloGoodEx}
Fix a group $H$. Recall from Subsection \ref{GraphWPSubsect} 
that the assignments: 
$$X \mapsto \bigoplus_X H \text{ and } X \mapsto \ast_X H $$
give rise to set-theoretic haloes of groups. Both satisfy (a) and (b). 
\end{Ex}

\begin{Lemma} \label{AutHaloLem}
    If $\mathcal{L}_X$ is a halo of groups over $X$ satisfying (a), 
then $\mathcal{AL}_X$ is a halo of groups over $X$. 
\end{Lemma}

\begin{proof}
    By Remark \ref{AutMostlyHalo}, it is sufficient to observe 
    that for $Y,Z \subseteq X$, every $\beta \in \AL_{f,X} (Y) \cap \AL_{f,X} (Z)$ 
    acts as the identity on both $L(X \setminus Y)$ and $L(X \setminus Z)$, 
    hence acts as the identity on 
    $\langle L(X \setminus Y) \cup L(X \setminus Z) \rangle = L (X \setminus (Y \cap Z))$. 
\end{proof}

Now suppose that $\mathcal{L}$ is a set-theoretic (respectively graph-theoretic) 
halo of groups satisfying (a) and (b). 
Let $X$ and $Y$ be sets (graphs) and let $\iota : Y \hookrightarrow X$ 
be an inclusion of $Y$ into $X$ as a subset (induced subgraph) of $X$. 
If we let $\beta \in \AL_f (Y) \leq \Aut (L(Y))$ then 
there is some $Z \subseteq Y$ finite for which $\beta$ restricts to 
an automorphism of $L(Z)$ and acts as the identity on $L(Y \setminus Z)$. 
Identifying $Y$ and $Z$ with their images in $X$, and applying condition (b), 
we obtain $\tilde{\beta} \in \AL_{f,X} (Z) \leq \AL_f (X)$ extending $\beta$. 
Set $\hat{\iota}(\beta) = \tilde{\beta} \in \AL_f (X)$ 
and note that $\hat{\iota} : \AL_f (Y) \rightarrow \AL_f (X)$ is a well-defined map 
by Remark \ref{AutOnePtGenRem}, and is clearly injective 
(since we may recover $\beta$ from $\tilde{\beta}$ by restricting to $Y$). 

\begin{Prop} \label{AutosAreHaloesProp}
If $\mathcal{L}$ is a set-theoretic 
(respectively graph-theoretic) halo of groups, 
satisfying (a) and (b) 
then the assignment $X \mapsto \AL_f (X)$ 
and $(\iota : Y \hookrightarrow X) \mapsto (\hat{\iota} : \AL_f (Y) \hookrightarrow \AL_f (X))$ 
(with $\hat{\iota}$ defined as above) 
gives rise to a set-theoretic 
(respectively graph-theoretic) halo of groups $\mathcal{AL}$. 
\end{Prop}

\begin{proof}
    By the discussion above and recalling Lemma \ref{AutHaloLem}, 
    it is sufficient to check that the above assignments yield a covariant functor, 
    which is immediate from definitions. 
\end{proof}

\begin{Prop} \label{AutoRFProp}
Let $Y \subseteq X$ and suppose $L(Y)$ is finitely generated residually finite. 
Then $\AL_{X,f} (Y)$ is residually finite. 
\end{Prop}

\begin{proof}
It is a classical result of Mal'cev that the automorphism 
group of any finitely generated residually finite group is residually finite. 
It is clear that restriction from $L(X)$ to $L(Y)$ induces an embedding 
of $\AL_{X,f} (Y)$ into $\Aut(L(Y))$, and the result follows. 
\end{proof}

We immediately deduce:

\begin{Cor} \label{AutEnrRFCoroll}
    Suppose that $\mathcal{L}$ is a set-theoretic 
    (respectively graph-theoretic) halo of groups satisfying conditions (a) and (b) above, 
    with the property that for every finite set (graph) $Z$, 
    $L(Z)$ is residually finite. 
    Then for every set (graph) $X$, $\AL_f(X)$ is locally residually finite, 
    hence LEF, hence sofic. 
\end{Cor}

\begin{Thm} \label{SoficAutEnrThm}
Suppose that $\Gamma$ is a sofic 
(respectively $\mathbb{K}$-linear sofic; hyperlinear; weakly sofic) group 
and that $\Gamma \curvearrowright_{\alpha} X$ is a sofic action 
on the set (respectively graph) $X$. 
    Suppose that $\mathcal{L}$ is a set-theoretic 
    (graph-theoretic) halo of groups satisfying conditions (a) and (b) above. 
    If $\AL_f(X)$ is a sofic 
(respectively $\mathbb{K}$-linear sofic; hyperlinear; weakly sofic) group, 
then so is $\mathcal{AL} (\Gamma,X)$. 
In particular, 
this holds if $L(Y)$ is finitely generated residually finite for every finite subset 
(induced subgraph) $Y \subseteq X$. 
\end{Thm}

\begin{proof}
    By Proposition \ref{AutosAreHaloesProp}, 
    $\mathcal{AL}$ is a set-theoretic (respectively graph-theoretic) halo of groups, 
    so Corollary \ref{BigHaloCorollSpec} is applicable. 
    The final statement follows from Corollary \ref{AutEnrRFCoroll}. 
\end{proof}

\begin{Thm} \label{LEFAutEnrThm}
Suppose that $\Gamma$ is a LEF group 
and that $\Gamma \curvearrowright_{\alpha} X$ is a LEF action 
on the set (respectively graph) $X$. 
    Suppose that $\mathcal{L}$ is a set-theoretic 
    (graph-theoretic) halo of groups satisfying conditions (a) and (b) above. 
    If $\AL_f(X)$ is a LEF group, 
then so is $\mathcal{AL} (\Gamma,X)$. 
In particular, 
this holds if $L(Y)$ is finitely generated residually finite for every finite subset 
(induced subgraph) $Y \subseteq X$. 
\end{Thm}

\begin{proof}
    This is much the same as the proof of Theorem \ref{SoficAutEnrThm}, 
    using Corollary \ref{LEFHaloProdCor} in place of Corollary \ref{BigHaloCorollSpec}. 
\end{proof}

\begin{Ex} \label{AutEnrApplEx}
suppose $\Gamma \curvearrowright_{\alpha} X$ is a sofic action (respectively LEF action) 
of the group $\Gamma$ on the set $X$. 
\begin{itemize}
    \item[(i)] If $H$ is a finitely generated residually finite group, 
    then for any finite set $Y$, $\ast_Y H$ 
    is finitely generated residually finite. 
    Moreover, as noted in Example \ref{SumProdHaloGoodEx}, 
    the assignment $Z \mapsto \ast_Z H$ is a set-theoretic halo of groups 
    satisfying assumptions (a) and (b). 
    Hence if $\Gamma \in \mathcal{C}$ for $\mathcal{C}$ 
    the class of sofic; $\mathbb{K}$-linear sofic; hyperlinear or weakly sofic groups 
    (respectively $\mathcal{C}$ the class of LEF groups) 
    then $\FAut (\ast_X H) \rtimes_{\hat{\alpha}} \Gamma \in \mathcal{C}$ also, 
    where $\FAut (\ast_X H)$ denotes the group of automorphisms of 
    $\ast_X H$ which are supported on $\ast_Y H \leq \ast_X H$ 
    for some finite subset $Y \subseteq X$. 
    In particular, taking $H = \mathbb{Z}$, we have 
    $\FAut (F(X)) \rtimes_{\hat{\alpha}} \Gamma \in \mathcal{C}$; 

    \item[(ii)] Similarly, if $H$ is a finitely generated residually finite group, 
    then so is $\bigoplus_Y H$, for any finite set $Y$. 
    Thus for $\mathcal{C}$ 
    the class of sofic; $\mathbb{K}$-linear sofic; hyperlinear or weakly sofic groups 
    (respectively $\mathcal{C}$ the class of LEF groups), 
    if $\Gamma \in \mathcal{C}$ then 
    $\FAut (\bigoplus_X H) \rtimes_{\hat{\alpha}} \Gamma \in \mathcal{C}$ also, 
    where $\FAut (\bigoplus_X H)$ denotes the group of automorphisms of 
    $\bigoplus_X H$ which are supported on $\bigoplus_Y H \leq \bigoplus_X H$ 
    for some finite subset $Y \subseteq X$. 

    \item[(iii)] As a special case of (ii), 
    if $R$ is a commutative unital ring whose underlying additive group $A$ is finitely generated, 
    then for any set $X$, $\GL_f (X,R)$ embeds as a subgroup of $\FAut (\bigoplus_X A)$, 
    and the action of $\Gamma$ on $\FAut (\bigoplus_X A)$ 
    restricts to the action on $\GL_f (X,R)$ described in Subsection \ref{LinEnrSubsect}. 
    Thus $\mathcal{GL}(\Gamma,X,R)$ 
    embeds as a subgroup of $\FAut (\bigoplus_X A) \rtimes_{\hat{\alpha}} \Gamma$. 
    Moreover $A$ is finitely generated residually finite. 
    Thus from (ii) we obtain sufficient hypotheses for 
    $\mathcal{GL}(\Gamma,X,R)$ to be 
    sofic; $\mathbb{K}$-linear sofic; hyperlinear; weakly sofic or LEF, 
    which are complementary to those of Theorems \ref{LinEnrSoficThm} and \ref{LinEnrLEFThm}. 
    
\end{itemize}
\end{Ex}



\subsection*{Acknowledgments} We would like to thank Francesco Fournier Facio and Srivatsav Kunnawalkam Elayavalli for enlightening conversations which informed our thinking 
on the subject of this paper. 

\bibliographystyle{alpha}
\bibliography{literature}

\end{document}